\documentclass[12pt,a4paper]{amsart}
\usepackage{textcomp}
\usepackage{indentfirst}
\usepackage{amsmath}
\usepackage{amsthm}
\usepackage{amsfonts}
\usepackage{amssymb}
\usepackage{mathrsfs}
\usepackage{latexsym}
\usepackage{amssymb}
\usepackage{amsfonts}
\usepackage{textcomp}
\usepackage{xcolor}
\usepackage{indentfirst}
\usepackage{mathrsfs}
\usepackage[all]{xy}
\usepackage{bbm}
\usepackage[T1]{fontenc}
\usepackage{enumerate}
\newtheorem{theorem}{Theorem}[section]
\newtheorem{lemma}[theorem]{Lemma}
\newtheorem{proposition}[theorem]{Proposition}
\newtheorem{corollary}[theorem]{Corollary}

\theoremstyle{definition}

\newtheorem{remark}[theorem]{Remark}

\newtheorem{example}[theorem]{Example}
\numberwithin{equation}{section}

\newcommand{\ta}{\theta,\alpha}
\newcommand{\at}{\alpha,\theta}
\newcommand{\bt}{\bar\theta}
\newcommand{\ba}{\bar\alpha}
\newcommand{\bv}{\bar\varphi}
\newcommand{\kaz}{k^{\alpha}_0}
\newcommand{\ktz}{k^{\theta}_0}
\newcommand{\bg}{\bar\gamma}

\def\ata{A^{\theta,\alpha}}
\def\aat{A^{\alpha,\theta}}
\def\nk0t{\|\tilde k_0^\theta\|^{-2}}
\def\kda{K_\alpha}
\def\kdt{K_\theta}

\def\b1{\mathcal{B}_1(\kdt,\kda)}

\def\lcm{\operatorname{lcm}}
\def\gcd{\operatorname{gcd}}

\newcommand{\dc}{P^- D_{|\theta H^2}}
\newcommand{\dd}{P^- D_{| H^2_-}}

\newcommand{\dtha}{P_{\alpha H^2}D_{|\theta H^2}}

\newcommand{\dba}{P_{\alpha H^2}D_{|H^2_-}}
\title[Intertwining property for compressions  ]{Intertwining property for  compressions of multiplication operators}

\author[M. C. C\^amara, K. Kli\'s--Garlicka, B. \L anucha, and M. Ptak]{M. Cristina C\^amara, Kamila Kli\'s--Garlicka, Bartosz \L anucha, and Marek Ptak}

\address{M. Cristina C\^{a}mara, Center for Mathematical Analysis, Geometry and Dynamical Systems\\ Mathematics
Department, Instituto Superior T\'{e}cnico, Universidade de Lisboa\\ Av. Rovisco Pais, 1049-
001 Lisboa, Portugal.}\email{ccamara@math.ist.utl.pt}

\address{Kamila Kli\'s--Garlicka, Department of Applied Mathematics,
University of Agriculture, ul. Balicka 253c\\ 30-198 Krak\'ow, Poland.}
\email{rmklis@cyfronet.pl}

\address{
Bartosz {\L}anucha, Institute of Mathematics,
 Maria Curie-Sk{\l}odowska University, pl. M.
Curie-Sk{\l}odowskiej 1, 20-031 Lublin, Poland}
\email{bartosz.lanucha@poczta.umcs.lublin.pl}

\address{Marek Ptak, Department of Applied Mathematics,
University of Agriculture, ul. Balicka 253c\\ 30-198 Krak\'ow, Poland.}\email{rmptak@cyf-kr.edu.pl}

\thanks{The work of the first author
was partially supported by FCT/Portugal through UID/MAT/04459/2020. The research of the second and the fourth authors was financed by the Ministry of Science and Higher Education of the Republic of Poland.}
\date{\today}

\begin{document}

\begin{abstract}
Following Beurling's theorem the natural compressions of the multiplication operator in the classical  $L^2$ space are compressions to model spaces and to their orthogonal complements. Two possibly different model spaces are considered hence asymmetric truncated Toeplitz and asymmetric dual truncated Toeplitz operators are investigated. The main purpose of the paper is to characterize operators which intertwine compressions of the unilateral shift.
\end{abstract}
\maketitle
\section{Introduction}%
Let $L^2$ be the classical space of (equivalent classes of) functions on the unit circle $\mathbb{T}$, which are measurable and square integrable with respect to the normalized Lebesgue measure $m$ on $\mathbb{T}$. Let $H^2$ be  the Hardy space, that is, the space of those functions from $L^2$, which can be extended to functions  holomorphic on whole unit disk $\mathbb{D}$. Denote by $P$ the natural projection from $L^2$ onto $H^2$. For any $\varphi\in L^\infty$ define a (classical) {\it Toeplitz operator} with symbol $\varphi$ as $T_\varphi f=P(\varphi f)$ for all $f\in H^2$. The unitary equivalence of the unilateral shift  $S$ and the Toeplitz operator $T_z$ makes Toeplitz operators especially interesting.

 Beurling's theorem characterizes all nontrivial invariant subspaces for $T_{\bar z}$ as $K_\theta:=H^2\ominus \theta H^2$ (these subspaces are called {\it model spaces}) with  a nonconstant inner function $\theta$ (for a constant inner function $\theta$ we will use the convention $\kdt=\{0\}$). Recently,  compressions of Toeplitz operators to model spaces $\kdt$ or their orthogonal complements $K_{\theta}^{\perp}=L^2\ominus \kdt$ were strongly investigated, see e.g. % \cite{sar1, CKLPd, CKLPcom, ding, hu, Ding Qin Sang}.
 \cite{sar1, BCT, ding, gar3, MYZ, hu, Ding Qin Sang}. In what follows $P_{\theta}$ will denote the orthogonal projection from $L^2$ onto
 $K_{\theta}$ and $P_{{\theta}}^\perp=I_{L^2}-P_{\theta}$ from $L^2$ onto $K_{\theta}^{\perp}$.

For  $\varphi\in L^2$  let $M_\varphi :D(M_\varphi)\to L^2$ be the densely defined multiplication operator $M_\varphi f=\varphi f$, where $D(M_\varphi)=\{f\in L^2\colon \varphi f\in L^2\}$. Note that $L^\infty\subset D(M_\varphi)$ for all $\varphi\in L^2$.
Recall after \cite{GMR} that $K_\theta^\infty:=K_\theta\cap L^\infty $ is a dense subset of $K_\theta$. %Since $\bar z\overline{H^\infty}$ is a dense subset of $H^2_-$ and $\theta H^\infty $ is a dense subset of $\theta H^2$, it follows that
The space $K_\theta^\perp\cap L^\infty $ is also a dense subset of $K_\theta^\perp$ as it was observed in \cite{CKLPcom}.
 %(For constant inner function $\theta $ we set $\kdt\{0\}$.)
 For nonconstant inner functions $\theta, \alpha$ and for $\varphi\in L^2$ define
$$A_{\varphi}^{\theta,\alpha}=P_{\alpha}M_{\varphi|K_{\theta}\cap L^{\infty}},\  B_{\varphi}^{\theta,\alpha}=P_{{\alpha}}^\perp M_{\varphi|K_{\theta}\cap L^{\infty}}\ \text {and}\  D_{\varphi}^{\theta,\alpha}=P_{{\alpha}}^\perp M_{\varphi|K_{\theta}^{\perp}\cap L^{\infty}}.$$
If $A_{\varphi}^{\theta,\alpha}$ extends to the whole $K_\theta$ as a bounded operator, it is called  an {\it asymmetric truncated Toeplitz operator} (ATTO, see \cite{part, part2, BCKP, blicharz1}). Similarly, if $B_{\varphi}^{\theta,\alpha}$ extends to a bounded operator from $K_\theta$ to $K_\alpha^\perp$,  it is called  a {\it big asymmetric truncated Hankel  operator} (ATHO, see \cite{MYZ, CKLPcom}), and if $D_{\varphi}^{\theta,\alpha}$ extends to the whole $K_\theta^\perp$ as a bounded operator, it is called  an {\it asymmetric dual truncated Toeplitz operator} (ADTTO, see \cite{CKLPd, CKLPcom}). It is easy to verify that $(A_{\varphi}^{\theta,\alpha})^*=A_{\bar\varphi}^{\alpha,\theta}$, $(D_{\varphi}^{\theta,\alpha})^*=D_{\bar\varphi}^{\alpha,\theta}$ and $(B_{\varphi}^{\theta,\alpha})^*=P_{\theta}M_{\bar\varphi|K_{\alpha}\cap L^{\infty}}$. Let us fix the notation
\begin{align*}
  \mathcal{T}(K_\theta,K_\alpha)&=\{A_{\varphi}^{\theta,\alpha}\colon\, \varphi\in
  L^2\ \mathrm{and}\ A_{\varphi}^{\theta,\alpha}\
  \mathrm{is\ bounded}\},\\
\mathcal{T}(K_\theta,K^\perp_\alpha)&=\{B_{\varphi}^{\theta}\colon\, \varphi\in
  L^2\ \mathrm{and}\ B_{\varphi}^{\theta,\alpha}\
  \mathrm{is\ bounded}\},\\ \mathcal{T}(K^\perp_\theta, K^\perp_\alpha)&=\{D_{\varphi}^{\theta,\alpha}\colon\, \varphi\in
  L^2\ \mathrm{and}\ D_{\varphi}^{\theta,\alpha}\
  \mathrm{is\ bounded}\}%=\{D_{\varphi}^{\theta}\colon\, \varphi\in L^\infty\}
  .\end{align*}
In case $\theta=\alpha$ we will use the shorter notation
$A_{\varphi}^{\theta}$, $ B_{\varphi}^{\theta}$, $D_{\varphi}^{\theta}$ and $\mathcal{T}(K_\theta)$ and  $\mathcal{T}(K_\theta,K^\perp_\theta)$, $\mathcal{T}(K^\perp_\theta)$, respectively.

Since the unilateral shift $S$ is unitarily equivalent to the Toeplitz operator $T_z$ we are able to describe the commutant of the unilateral shift as
\[\{S\}^\prime=\{T_z\}^\prime=\{T_\varphi: \varphi\in H^\infty:=L^\infty\cap H^2\}.\]

Now considering the compressions $A^\theta_z$ and $D^\theta_z$, it is natural to try to describe the commutants $\{A^\theta_z\}^\prime$ and   $\{D^\theta_z\}^\prime$. In the more general asymmetric setting we are searching for all operators intertwining  $A^\theta_z$ and  $A^\alpha_z$ in the case of two model spaces or searching for all operators intertwining  $D^\theta_z$ and  $D^\alpha_z$ in the case of orthogonal complements of two model spaces, i.e., we try to describe the following sets of operators
\[\mathcal{I}(\kdt,\kda)=\{B\in \mathcal{B}(\kdt,\kda)\colon A^\alpha_z B=B A^\theta_z\},\]
\[\mathcal{I}(\kdt^\perp,\kda^\perp)=\{B\in \mathcal{B}(\kdt^\perp,\kda^\perp)\colon D^\alpha_z B=B D^\theta_z\}.\]

Section 3 is devoted to characterization of operators intertwining  $A^\alpha_z$ and $A^\theta_z$. Such characterization (but less precise) was given in \cite{Bercovici} (Theorem III.1.16). However, the proof there is based on the commutant lifting theorem. Here we give a more basic proof.

In Section 4 we investigate functional calculus for asymmetric truncated Toeplitz operators.

In Section 5, Theorem \ref{idatto}, we characterize those asymmetric dual truncated Toeplitz operators which intertwine $D^\alpha_z$ and $D^\theta_z$.
In Section~5 we give a necessary and sufficient condition for an operator from $\mathcal{B}(K_\theta^\perp, K_\alpha^\perp)$ to intertwine  $D^\alpha_z$ and $D^\theta_z$ (Theorem \ref{kmutant}). Our results from Sections 4 and 5 were obtained as new even in the symmetric case $\theta=\alpha$. However, this special case was considered in the very recent papers \cite{Gu2,LSD}.

\section{Operators with analytic symbols}

Recall the standard notations  for a nonconstant inner function $\theta$. Namely $k_0^\theta=P_\theta 1=1-\overline{\theta(0)}\theta$ .
We also set $S_\theta=A_z^\theta$.

\begin{proposition}\label{pr21}
Let $\alpha$, $\theta$ be nonconstant inner functions and let $\varphi\in H^\infty$.
\begin{enumerate}
\item If $\alpha\leqslant \theta$ (i.e., $\alpha$ divides $\theta$), then $(\ata_\varphi)^*=\aat_{\bar\varphi}=T_{\bar\varphi |\kda}$.
\item $\ata_\varphi=0$ if and only if $\varphi\in \alpha H^\infty$.
\end{enumerate}
\end{proposition}

\begin{proof}
It is easy to see that $(A_{\varphi}^{\theta,\alpha})^*=A_{\bar\varphi}^{\alpha,\theta}=P_{\theta}T_{\bar\varphi |\kda}$. Since $\kda\subset K_{\theta}$ is invariant for $T_{\bar\varphi}$, we get $P_\theta T_{\bar\varphi|\kda}=T_{\bar\varphi|\kda}$. Hence (1) follows.

Part (2) follows from \cite[Thm 2.1]{blicharz1}.
\end{proof}

\begin{corollary}\label{cor1}
Let $\alpha$, $\theta$ be nonconstant inner functions and $\varphi_1,\varphi_2\in H^\infty$. Then
$$\ata_{\varphi_1}=\ata_{\varphi_2} \text{ if and only if } \varphi_1-\varphi_2\in\alpha H^\infty.$$
\end{corollary}

Let $\mathcal{A}(\theta,\alpha)$ denote the set of all bounded asymmetric truncated Toeplitz operators in $\mathcal{B}(K_\theta,K_\alpha)$ with analytic symbols, i.e.,
$$\mathcal{A}(\theta,\alpha)=\{A^{\theta,\alpha}_\varphi\in \mathcal{T}(K_\theta,K_\alpha)\colon \varphi\in H^2\}.$$
\begin{remark}\label{mk}
If  $A^{\theta,\alpha}_\varphi\in\mathcal{A}(\theta,\alpha)$, then by \cite[Corollary 2.6]{blicharz1} there is $\psi\in K_\alpha$ such that  $A^{\theta,\alpha}_\varphi=A^{\theta,\alpha}_\psi$.
\end{remark}

% Denote $S_\theta=A_z^\theta$ and $S_\alpha=A_z^\alpha$.
\begin{lemma}
Let $A\in \mathcal{B}(K_\theta,K_\alpha)$. Then the following are equivalent
\begin{enumerate}
\item $A\in \mathcal{A}(\theta,\alpha)$;
\item there is $\psi\in K_\alpha$ such that $A-S_\alpha A S^*_\theta=\psi\otimes k_0^\theta$;
\item $\langle AS^* f, S^* g\rangle=\langle Af,g\rangle$ for all $f\in K_\theta$, $g\in K_\alpha$ such that $f(0)=0$.
\end{enumerate}
\end{lemma}

\begin{proof}
It follows from \cite[Theorem 2.1]{BM} and its proof that if $A$ satisfies (2), then $A=A^{\theta,\alpha}_\psi\in \mathcal{A}(\theta,\alpha)$. On the other hand, if $A\in \mathcal{A}(\theta,\alpha)$, then, by Remark \ref{mk},  we have $A=A^{\theta,\alpha}_\psi$ with $\psi\in K_\alpha$. By \cite[Lemma 2.2]{BM},
$$A-S_\alpha A S^*_\theta=(\psi-\psi(0)k_0^\alpha)\otimes k_0^\theta+k_0^\alpha \otimes (\overline{\psi(0)} k_0^\theta)=\psi\otimes k_0^\theta.$$
Thus (1) and (2) are equivalent.

Assume now that $$A-S_\alpha A S^*_\theta=\psi\otimes k_0^\theta.$$ Then for $f\in K_\theta$, $g\in K_\alpha$ such that $f(0)=0$ it holds
\begin{multline*}
\langle Af,g\rangle-\langle AS^* f, S^* g\rangle
=\langle Af,g\rangle-\langle AS^*_\theta f, S^*_\alpha g\rangle=\langle(A-S_\alpha A S^*_\theta)f,g \rangle\\= \langle (\psi\otimes k_0^\theta)f,g\rangle=\langle f, k_0^\theta\rangle \langle\psi,g\rangle=f(0) \langle\psi,g\rangle=0.\end{multline*}
Thus $(2)\Rightarrow(3)$.

To prove the implication $(3)\Rightarrow (2)$ define $B=A-S_\alpha A S^*_\theta$. Then by (3), for any $c\in \mathbb{C}$ we have
$$0=B(I_{K_\theta}-ck_0^\theta\otimes k_0^\theta)=B-(cBk_0^\theta)\otimes k_0^\theta.$$
and (2) follows since $cBk_0^\theta\in K_\alpha$.
\end{proof}

\section{Intertwining with restricted shifts}
Let $\alpha$, $\theta$ be nonconstant inner functions. Denote by $\lcm(\alpha, \theta)$ and by $\gcd(\alpha, \theta)$ the least common multiplier and the greatest common divisor of $\theta$ and $\alpha$, respectively.

Recall that $\mathcal{I}(K_\theta,K_\alpha)$ is the set of all operators  $A\in \mathcal{B}(K_\theta,K_\alpha)$ which intertwine $S_\alpha$ and $S_\theta$, i.e.,
\begin{equation}\label{relation}
S_\alpha A=A S_\theta.
\end{equation}
\begin{remark}
	Note that the operators satisfying \eqref{relation} were studied and characterized in \cite{Bercovici} (see Theorem III.1.16 and Proposition III.1.21 in \cite{Bercovici}). However, the proof there is based on the commutant lifting theorem. Below we obtain a characterization of these operators using more basic methods (Theorem \ref{inter}).
\end{remark}
To describe  $\mathcal{I}(K_\theta,K_\alpha)$ we start with the following lemma.

\begin{lemma}\label{3.11}
Let $\alpha, \theta$ be  nonconstant inner functions. Then
\begin{equation}\alpha K_\theta \cap \theta K_\alpha=\lcm(\alpha, \theta)\, K_{\gcd(\alpha, \theta)}.
\end{equation}
\end{lemma}

\begin{proof}
Let $\gamma=\gcd(\alpha,\theta)$, $\eta=\lcm(\alpha,\theta)$ and put $\alpha_1=\frac{\alpha}{\gamma}$, $\theta_1=\frac{\theta}{\gamma}$. Since $\alpha_1$, $\theta_1$ are relatively prime (have no nonconstant common divisor), we have
\begin{align*}
\alpha K_\theta \cap \theta K_\alpha = & \alpha (K_{\theta_1}\oplus \theta_1 K_\gamma) \cap \theta (K_{\alpha_1}\oplus \alpha_1 K_\gamma) =\\
=& (\alpha K_{\theta_1}\oplus \tfrac{\alpha\theta}{\gamma} K_\gamma)\cap (\theta K_{\alpha_1}\oplus \tfrac{\alpha\theta}{\gamma} K_\gamma)\\
=&\gamma (\alpha_1 K_{\theta_1}\cap \theta_1 K_{\alpha_1})\oplus \eta K_\gamma.
\end{align*}
Let $f\in \alpha_1 K_{\theta_1}\cap
\theta_1K_{\alpha_1}$. Then $f=\alpha_1f_1=\theta_1f_2$, where $f_1\in K_{\theta_1}$, $f_2\in K_{\alpha_1}$. Since $\alpha_1,\theta_1$ are relatively prime, then $f_1=\theta_1f_1^\prime$ with $f_1^\prime\in H^2$. This is possible only if $f=f_1=0$, which means that $\alpha_1 K_{\theta_1}\cap \theta_1K_{\alpha_1}=\{0\}$.
\end{proof}
\begin{example}

Take $a\in\mathbb{D}\setminus\{0\}$ and let $\theta(z)=z\varphi_{a}(z)=z\frac{a-z}{1-\bar a z}$, $\alpha(z)=z^2$. Then $\gamma(z)=z$ and
$$\alpha K_\theta \cap \theta K_\alpha=z^2\varphi_{a}K_z=z^2\varphi_{a}\mathbb{C}.$$
\end{example}
\begin{corollary}
Let $\alpha$, $\theta$ be  nonconstant inner functions such that $\alpha \leqslant \theta$. Then
$\alpha K_\theta \cap \theta K_\alpha= \theta\, K_{\alpha}.
$%
\end{corollary}

\begin{corollary}
Let $\alpha$, $\theta$ be nonconstant relatively prime inner functions. Then $\alpha K_\theta \cap \theta K_\alpha= \{0\}.
$
\end{corollary}
\bigskip

Observe that as in \cite{sar1}, if $A\in\mathcal{I}(K_\theta,K_\alpha)$, then
\begin{equation}
\label{warr} A-S_\alpha A S_\theta^*=A(I_{K_{\theta}}- S_\theta S_\theta^*)=(Ak_0^{\theta})\otimes k_0^{\theta},
\end{equation}
and by \cite[Theorem 2.1]{BM}, $A\in\mathcal{T}(K_\theta,K_\alpha)$. % and $A= A_{\varphi}^{\theta,\alpha}$ with $\varphi=Ak_0^{\theta}\in K_{\alpha}$.
Therefore we have that $\mathcal{I}(K_\theta,K_\alpha)\subset\mathcal{T}(K_\theta,K_\alpha)$. To determine which  $A_{\varphi}^{\theta,\alpha}$ belong to $\mathcal{I}(K_\theta,K_\alpha)$ we will study the difference
\begin{equation*}%\label{3.1}
  S_\alpha A_{\varphi}^{\theta,\alpha}-A_{\varphi}^{\theta,\alpha} S_\theta.
\end{equation*}
Recall firstly that for any inner $\theta$ the  operator $C_{\theta}:L^2\to L^2$ defined by
$$C_{\theta}f(z)=\theta\bar z\overline{f(z)}$$
is an antilinear isometry such that $C_{\theta}^2=I_{L^2}$. Moreover, $C_{\theta}$ preserves $K_{\theta}$ and maps $\theta H^2$ onto $\bar z\overline{H^2}$ and vice versa \cite{GMR, GP}. We will use the notation $\tilde{k}_0^\theta$ for $C_{\theta}{k}_0^\theta$.

Let $f\in K_\theta^{\infty}$ and $g\in K_\alpha^{\infty}$. Since $P_{\alpha}=P-M_{\alpha} PM_{\bar\alpha}$ on $L^2$ the following holds
\begin{align*}
\langle S_\alpha A_{\varphi}^{\theta,\alpha}f,g \rangle &= \langle S_\alpha P_\alpha(\varphi f),g\rangle = \langle z P_\alpha(\varphi f),g\rangle  \\
&=\langle zP(\varphi f),g\rangle -\langle z\alpha P_\alpha(\bar\alpha\varphi f),g\rangle\\
&=\langle \varphi f,P(\bar zg)\rangle=\langle \varphi f,\bar zg-\bar zg(0)\rangle\\
&=\langle \varphi z f,g\rangle-\langle f, P_{\theta}(\bar z \bar \varphi)\rangle  \langle {k}_0^\alpha, g\rangle.
\end{align*}
On the other hand, since $S_\theta f=zf-\langle f,\tilde{k}_0^\theta\rangle\theta$, we have
$$\langle A_{\varphi}^{\theta,\alpha} S_\theta f,g \rangle = \langle \varphi zf,g\rangle - \langle f,\tilde{k}_0^\theta \rangle \langle \varphi\theta, g\rangle. $$
Therefore we have the following relation
\begin{equation}\label{3.2}
S_\alpha A_{\varphi}^{\theta,\alpha}-A_{\varphi}^{\theta,\alpha} S_\theta=P_{\alpha}(\varphi \theta)\otimes \tilde{k}_0^\theta -{k}_0^\alpha\otimes P_{\theta}(\bar z \bar \varphi)
\end{equation}
(compare with \cite[Corollary 2.4(d)]{BM}).

\begin{lemma}\label{3.3}The right hand side of \eqref{3.2} cancels, i.e.,
\begin{equation}\label{3.4}P_{\alpha}(\varphi \theta)\otimes \tilde{k}_0^\theta ={k}_0^\alpha\otimes P_{\theta}(\bar z \bar \varphi)\end{equation}
if and only if there is $c\in \mathbb{C}$ such that $(\varphi\theta -c)\perp K_{\lcm(\alpha,\theta)}$.
\end{lemma}

\begin{proof} Note that \eqref{3.4} holds if and only if there is  $c\in \mathbb{C}$ such that
\[P_{\alpha}(\varphi \theta)=c\,{k}_0^\alpha\quad\text{and}\quad P_{\theta}(\bar z \bar \varphi)= \bar c\,\tilde{k}_0^\theta.\]
The first condition is equivalent to $P_{\alpha}(\varphi \theta-c)=0$, which means that $(\varphi \theta-c)\, \perp K_\alpha$.
The second condition is equivalent to
\[P_{\theta}( \varphi\theta)=  P_{\theta}(C_\theta( \bar z\bar \varphi))=C_\theta P_{\theta}(\bar z \bar \varphi)=
C_\theta (\bar c\,\tilde{k}_0^\theta)= c\, k_0^\theta, \]
which, similarly as above, gives that $(\varphi\theta - c)\,\perp K_{\theta}$. Hence \eqref{3.4} holds if and only if $$(\varphi\theta - c)\,\perp \operatorname{span}\{K_\alpha, K_\theta\}=K_{\lcm(\alpha,\theta)}$$ (see \cite[Corollary 5.9]{GMR}).% The reverse implication follows from the above.
\end{proof}

Observe that if $\varphi\in L^2$, then $(\varphi\theta -c)\perp K_{\lcm(\alpha,\theta)}$ (for some $c\in\mathbb{C}$) if and only if
$$(\varphi\theta -c)\in\overline{z H^2}+\lcm(\alpha,\theta)H^2,$$
or, equivalently,
$$\varphi\in \bar{\theta}\mathbb{C}+\overline{\theta z H^2}+\bar{\theta}\lcm(\alpha,\theta)H^2.$$
Note that the set on the right hand side is equal to
\begin{align*}
\overline{\theta H^2}+\bar{\theta}\lcm(\alpha,\theta)H^2&=\overline{\theta H^2}+\tfrac{\alpha}{\gcd(\alpha,\theta)}H^2\\
&=\overline{\theta H^2}+\alpha H^2+\tfrac{\alpha}{\gcd(\alpha,\theta)}K_{\gcd(\alpha,\theta)}.
\end{align*}

\begin{corollary}\label{wnnn}
	Let $\varphi\in L^2$. Then
	\begin{equation*}
	S_\alpha A_{\varphi}^{\theta,\alpha}-A_{\varphi}^{\theta,\alpha} S_\theta=0
	\end{equation*}
	if and only if
	$\varphi\in \overline{\theta H^2}+\alpha H^2+\tfrac{\alpha}{\gcd(\alpha,\theta)}K_{\gcd(\alpha,\theta)}$.
\end{corollary}

\begin{corollary}
	Every operator $A\in \mathcal{I}(K_{\theta}, K_{\alpha})$ has a symbol from $\tfrac{\alpha}{\gcd(\alpha,\theta)}K_{\gcd(\alpha,\theta)}$.
\end{corollary}

%{The intertwining relation \eqref{relation} was considered in \cite[Theorem III.1.16]{Bercovici}.  The proof there is based on the commutant lifting theorem. Here we give a more basic proof.}
\begin{theorem}\label{inter}
Let $\alpha, \theta$ be  nonconstant inner functions and let $A\in \mathcal{B}(K_{\theta}, K_{\alpha})$. Then $ S_\alpha A=A S_\theta$ if and only if $A\in\mathcal{T}(K_{\theta}, K_{\alpha})$ and $A=A_{\varphi}^{\theta,\alpha}$ with  $\varphi\in\tfrac{\alpha}{\gcd(\alpha, \theta)}\, K_{\gcd(\alpha, \theta)}$.
 %$$K_{\alpha}\cap (\alpha \overline{\theta}K_{\theta})=K_{\alpha}\cap (\alpha \overline{zK_{\theta}}).$$
\end{theorem}
\begin{remark}\label{3.10} The consequence of \cite[Theorem III.1.16]{Bercovici} is that there is $\varphi\in  \tfrac{\alpha}{\gcd(\alpha, \theta)}\, H^\infty$. Here we will show that $\varphi$ is in the smaller space but we can not guarantee that it is bounded. Moreover in the proof below we will not use the commutant lifting theorem.
\end{remark}
\begin{proof}
	If $ S_\alpha A=A S_\theta$, then by \eqref{warr} and \cite[Theorem 2.1]{BM}, $A\in\mathcal{T}(K_\theta,K_\alpha)$ and $A= A_{\varphi}^{\theta,\alpha}$ with $\varphi=Ak_0^{\theta}\in K_{\alpha}$. By Corollary \ref{wnnn}, $\varphi\in \overline{\theta H^2}+\alpha H^2+\tfrac{\alpha}{\gcd(\alpha,\theta)}K_{\gcd(\alpha,\theta)}$ and since $\theta$ is nonconstant, it follows that
	$$\varphi \in K_{\alpha}\cap\big(\overline{\theta H^2}+\alpha H^2+\tfrac{\alpha}{\gcd(\alpha,\theta)}K_{\gcd(\alpha,\theta)}\big)= \tfrac{\alpha}{\gcd(\alpha,\theta)}K_{\gcd(\alpha,\theta)}.$$

	On the other hand, if $A= A_{\varphi}^{\theta,\alpha}\in\mathcal{T}(\theta,\alpha)$ with $\varphi \in \tfrac{\alpha}{\gcd(\alpha, \theta)}\, K_{\gcd(\alpha, \theta)}$, then $ S_\alpha A=A S_\theta$ by Corollary \ref{wnnn}.
	\end{proof}

\begin{corollary}\label{inter2}
Let $\alpha, \theta$ be  nonconstant inner functions.
	Let $A\in \mathcal{B}(K_{\theta}, K_{\alpha})$. Then $ S_\alpha A=A S_\theta$ if and only if $A\in\mathcal{T}(\theta,\alpha)$ and
\begin{enumerate}
\item if $\alpha \leqslant \theta$, then $A=A_{\varphi}^{\theta,\alpha}$ has a symbol from  $ K_{\alpha}$,
\item if $\alpha \geqslant \theta$, then $A=A_{\varphi}^{\theta,\alpha}$ has a symbol from  $ \tfrac{\alpha}{\theta}K_{\theta}$,
    \item if $\alpha = \theta$, then $A=A_{\varphi}^{\theta}$ has a symbol from  $ K_{\theta}$.
\end{enumerate}
\end{corollary}%

%{Recall that (3) was shown in  \cite[Proposition III.1.21]{Bercovici}.}

\begin{corollary}
Let $\alpha, \theta$ be  nonconstant relatively prime inner functions. Then the only operator  $A\in \mathcal{B}(K_{\theta}, K_{\alpha})$ fulfilling  $ S_\alpha A=A S_\theta$ is the zero operator.
\end{corollary}

\begin{corollary}
	Let $\alpha, \theta$ be  nonconstant inner functions and let $A\in \mathcal{B}(K_{\theta}, K_{\alpha})$. Then $ S_\alpha^* A=A S_\theta^*$ if and only if $A\in\mathcal{T}(K_{\theta}, K_{\alpha})$ and $A=A_{\overline{ \varphi}}^{\theta,\alpha}$ with  $\varphi\in\tfrac{\alpha}{\gcd(\alpha, \theta)}\, K_{\gcd(\alpha, \theta)}$.
	%$$K_{\alpha}\cap (\alpha \overline{\theta}K_{\theta})=K_{\alpha}\cap (\alpha \overline{zK_{\theta}}).$$
\end{corollary}

In his unpublished paper \cite{Gu}, C. Gu investigated operators on $K_{\theta}$ intertwining $S_\theta$ and $S_\theta^*$. Let us observe here that for $B\in \mathcal{B}(K_{\theta}, K_{\alpha})$ the equation
\begin{equation}\label{STS}
S_\alpha B-BS_\theta^*=0
\end{equation}
can be equivalently expressed as
\begin{equation*}%\label{STS}
S_{\alpha^{\#}} A-AS_\theta^*=0,
\end{equation*}
where $A=J^{\#}BC_{\theta}\in \mathcal{B}(K_{\theta}, K_{\alpha^{\#}})$ and $J^{\#}:L^2\to L^2$,
$$J^{\#}f(z)=f^{\#}(z)=\overline{f(\overline{z})}.$$
Indeed, this follows easily from the fact that $J^{\#}$ is a conjugation which maps $K_{\alpha}$ onto $K_{\alpha^{\#}}$, $J^{\#}S_{\alpha^{\#}}J^{\#}=S_{\alpha}$ (see \cite{GMR} or \cite{BM}) and $S_\theta$ is $C_\theta$--symmetric, i.e., $S_\theta^*=C_\theta S_\theta C_\theta$.

Now, by Theorem \ref{inter}, $B\in \mathcal{B}(K_{\theta}, K_{\alpha})$ satisfies \eqref{STS} if and only if $J^{\#}BC_{\theta}\in \mathcal{T}(K_{\theta}, K_{\alpha^{\#}})$, $J^{\#}BC_{\theta}=A_{\varphi}^{\theta,\alpha^{\#}}$ with
$\varphi\in\tfrac{\alpha^{\#}}{\gcd(\alpha^{\#}, \theta)}\, K_{\gcd(\alpha^{\#}, \theta)}$. In other words, $B$ satisfies \eqref{STS} if and only if
$$B=J^{\#}A_{\varphi}^{\theta,\alpha^{\#}}C_\theta=P_{\alpha}J(I-P)M_{\psi|K_{\theta}}$$
(here $Jf(z)=\bar z f(\bar z)$) with
$$\psi=\overline{ \theta\varphi}\in\overline{\lcm(\alpha^{\#}, \theta)\cdot K_{\gcd(\alpha^{\#}, \theta)}}$$
(see \cite[Proposition 3.2(f)]{BM}). In \cite{Gu}, the operators of such form as above were called truncated Hankel operators (see also \cite{BM} and \cite{LM}). The result thus obtained is a generalization of Corollary 5.2(a) from \cite{Gu}.
\begin{corollary}
	Let $\alpha, \theta$ be  nonconstant inner functions and let $B\in \mathcal{B}(K_{\theta}, K_{\alpha})$. Then $ S_\alpha B=B S_\theta^*$ if and only if $B=P_{\alpha}J(I-P)M_{\psi|K_{\theta}}$ with $\psi\in\overline{\lcm(\alpha^{\#}, \theta)\cdot K_{\gcd(\alpha^{\#}, \theta)}}$.
\end{corollary}
It is also not difficult to verify that $B\in \mathcal{B}(K_{\theta}, K_{\alpha})$ satisfies
\begin{equation}\label{SST}
S_\alpha^* B-BS_\theta=0
\end{equation}
if and only if $C_\alpha BC_\theta$ satisfies \eqref{STS}, that is $B=C_{\alpha}J^{\#}A_{\varphi}^{\theta,\alpha^{\#}}$ with $\varphi\in\tfrac{\alpha^{\#}}{\gcd(\alpha^{\#}, \theta)}\, K_{\gcd(\alpha^{\#}, \theta)}$. Using \cite[Proposition 3.2(b)]{BM} we obtain
\begin{corollary}
	Let $\alpha, \theta$ be  nonconstant inner functions and let $B\in \mathcal{B}(K_{\theta}, K_{\alpha})$. Then $ S_\alpha^* B=B S_\theta$ if and only if $B=P_{\alpha}J(I-P)M_{\psi|K_{\theta}}$ with $\psi\in\overline{\gcd(\alpha^{\#}, \theta)}K_{\gcd(\alpha^{\#}, \theta)}=\bar z\overline{K_{\gcd(\alpha^{\#}, \theta)}}$.
\end{corollary}

\section{Functional calculus for ATTO's}

Since $A_z^{\theta,\alpha}$ is a completely non-unitary contraction we can refer to B. Sz.-Nagy--Foias functional calculus \cite[Theorem III.2.1]{NF}.

The weak* topology (ultraweak topology) in $\mathcal{B}(\kdt,\kda)$ is given by trace class operators of the form $\displaystyle t=\sum_{n=0}^\infty f_n\otimes g_n$ with $f_n\in\kdt$, $g_n\in\kda$ such that $\displaystyle\sum_{n=0}^\infty \|f_n\|^2< \infty$, $\displaystyle\sum_{n=0}^\infty \|g_n\|^2< \infty$.   Denote by $\mathcal{B}_1(\kdt,\kda)$ the space of all such trace class operators. Note that if $\displaystyle\sum_{n=0}^\infty f_n\otimes g_n\in\b1$, then $\displaystyle\sum_{n=0}^\infty f_n\bar g_n\in L^1$.

\begin{proposition}\label{p71} Let $\alpha$, $\theta$ be nonconstant inner functions. The mapping
 $\Phi\colon L^\infty\to \mathcal{T}(K_\theta,K_\alpha)$, where $\Phi(\varphi)=\ata_\varphi$, % $\varphi\mapsto\ata_\varphi$
 is a linear contractive  weak* -- weak*(WOT) continuous homomorphism.
\end{proposition}

%Define $\Phi\colon H^\infty\to \mathcal{A}(K_\theta, K_\alpha)$, where $\varphi\mapsto\ata_\varphi$.  We also  that $\Phi$ is weak* continues.
\begin{proof}The linearity is trivial. The contractivity follows from the inequalities:
\begin{equation}
\| \Phi(\varphi)\|=\| P_\alpha {T_\varphi}_{|{K_\theta}}\|\leqslant \| P_\alpha T_\varphi\|\leqslant\|T_\varphi\|\leqslant\|\varphi\|_\infty.
\end{equation}

To show continuity of $\Phi$ let us take $u_\beta\to u$ in the weak* topology in $L^\infty$. Then for any $t=\displaystyle\sum_{n=0}^\infty f_n\otimes g_n\in \mathcal{B}_1(\kdt,\kda)$ we have
$$\langle\Phi(u_\beta),t\rangle=\int_{\mathbb{T}} u_\beta\sum_{n=0}^\infty f_n\bar g_n\,dm\to \int_{\mathbb{T}} u\sum_{n=0}^\infty f_n\bar g_n\,dm=\langle \Phi(u),t\rangle.$$
Continuity in WOT topology can be proved similarly.
\end{proof}

 Recall that weak* topology on $ H^\infty\subset L^\infty$ is given by the duality of $L^\infty$ to $L^1$ ($(L^1)^*=L^\infty$).
 Denote by $H_0^1$ the subspace of $H^1\subset L^1$ containing functions vanishing at zero. Recall that the dual to $H^\infty$ is isomorphic to $L^1\Big/H^1_0$ since $H^\infty=(H^1_0)^\perp$. %  To consider the weak* continuity recall after \cite[p. 41]{Bercovici}
 Note that the following subspaces $\alpha H^\infty\subset \tfrac{\alpha}{\gcd(\alpha, \theta)}\, H^\infty\subset H^\infty\subset L^\infty$ are weak* closed and \[\alpha H^\infty =(\bar \alpha H^1_0)^\perp,\qquad \tfrac{\alpha}{\gcd(\alpha, \theta)}\, H^\infty=\Big(\overline{\tfrac{\alpha}{\gcd(\alpha, \theta)}}\, H^1_0\Big)^\perp.\] As a consequence of Hahn--Banach theorem we get
 \[\Big(\bar\alpha H^1_0\Big/\overline{\tfrac{\alpha}{\gcd(\alpha, \theta)}}\, H^1_0\Big)^*= \tfrac{\alpha}{\gcd(\alpha, \theta)}\, H^\infty\Big/\alpha H^\infty  \] and the equality is up to isometric isomorphism.
  The duality is given by the bilinear form
\begin{multline}\label{3.3a}\langle \,\tfrac{\alpha}{\gcd(\alpha, \theta)}\varphi+\alpha H^\infty,\bar\alpha f+\overline{\tfrac{\alpha}{\gcd(\alpha, \theta)}}\, H_0^1\,\rangle=\int_{\mathbb{T}} \tfrac{\alpha}{\gcd(\alpha, \theta)}\varphi \bar \alpha f\, dm\\=\int_{\mathbb{T}} \varphi  f\,\overline{\gcd(\alpha, \theta)}\, dm, \quad \varphi\in H^\infty, \ f\in  H_0^1.\end{multline}

 On the other hand, the space $\tfrac{\alpha}{\gcd(\alpha, \theta)}\, H^\infty\Big/{{\alpha H^\infty}}$ is isomorphic to  $  \tfrac{\alpha}{\gcd(\alpha, \theta)} \Big(H^\infty\Big/{\gcd(\alpha, \theta)}H^\infty\Big)$ and also $\bar\alpha H^1_0\Big/\overline{\tfrac{\alpha}{\gcd(\alpha, \theta)}}\, H^1_0$ is isomorphic to $\overline{\tfrac{\alpha}{\gcd(\alpha, \theta)}}\Big({{\gcd(\alpha, \theta)}}\, H^1_0\Big/ H^1_0\Big)$.
 Similarly as above,
 \[\Big(\overline{{\gcd(\alpha, \theta)}}\, H^1_0\Big/{ H^1_0}\Big)^*=  H^\infty\Big/{\gcd(\alpha, \theta)}H^\infty . \]
 Hence
 \begin{equation}\label{dual}\Big(\overline{\tfrac{\alpha}{\gcd(\alpha, \theta)}}\,\Big(\overline{{\gcd(\alpha, \theta)}} H^1_0\Big/ H^1_0\Big)\Big)^*= \tfrac{\alpha}{\gcd(\alpha, \theta)} \Big(H^\infty\Big/{\gcd(\alpha, \theta)}H^\infty \Big). \end{equation}
  The duality here is given by the bilinear form
\begin{multline}\langle\, \tfrac{\alpha}{\gcd(\alpha, \theta)}\big(\varphi+\gcd(\alpha,\theta) H^\infty\big), \overline{\tfrac{\alpha}{\gcd(\alpha, \theta)}}\big(\overline{\gcd(\alpha, \theta)} f+ H_0^1\big)\,\rangle\\=\int_{\mathbb{T}} \varphi  f\,\overline{\gcd(\alpha, \theta)}\, dm, \quad \varphi\in H^\infty, \ f\in  H_0^1.\end{multline}
Hence it coincides with \eqref{3.3a}.

\begin{theorem}\label{th72}
Let $\alpha$, $\theta$ be nonconstant inner functions  and let $X\in\mathcal{I}(K_\theta, K_\alpha)$.  Then there is $\varphi\in \tfrac{\alpha}{\gcd(\alpha, \theta)}\, H^\infty$ such that $X=\ata_\varphi$. Moreover, we have
\begin{equation}\label{eq210}
\|X\|=\operatorname{dist}(\varphi,\alpha H^\infty)\leqslant \|\varphi\|_\infty.
\end{equation}
Let $\Phi$ be as in Proposition \ref{p71}, $\pi$ be the quotient map $\pi\colon \frac \alpha{\gcd(\alpha,\theta)} H^\infty\to \frac \alpha{\gcd(\alpha,\theta)} H^\infty \Big/ \alpha H^\infty$ and $\iota_1$, $\iota_2$ be the natural inclusions. Then
there is an isometry \[\hat{\Phi}\colon \tfrac{\alpha}{\gcd(\alpha,\theta)}H^\infty/\alpha H^\infty \to \mathcal{I}(K_\theta, K_\alpha),\quad [\varphi]\mapsto\ata_\varphi,\]such that the following diagram is commuting, i.e., $\iota_2\circ\hat\Phi\circ \pi=\Phi\circ \iota_1$.

 \[\xymatrixcolsep{1pc}\xymatrix{
\frac \alpha{\gcd(\alpha,\theta)} H^\infty \ar[rr]^\pi \ar@{^{(}->}[rd]^{\iota_1} & & \frac \alpha{\gcd(\alpha,\theta)} H^\infty \Big/ \alpha H^\infty \ar[rr]^{\hat\Phi} & & \mathcal{I}(K_\theta, K_\alpha) \ar@{^{(}->}[dl]_{\iota_2}\\
& L^\infty \ar[rr]^\Phi & & \mathcal{T}(K_\theta, K_\alpha)
}\]

 The map $\hat\Phi$ is continuous with respect to the weak* topology in both spaces. % Moreover, $\mathcal{I}(K_\theta, K_\alpha)$ has property $\mathbb{A}_1(1)$.
\end{theorem}

\begin{proof}
The existence of $\varphi\in \frac{\alpha}{\gcd(\alpha,\theta)}H^\infty$ for a given $X\in \mathcal{I}(K_\theta, K_\alpha)$ such that $X=\ata_\varphi$ follows from \cite[Theorem III.1.16]{Bercovici}, see also  Remark \ref{3.10}. The inequalities in \eqref{eq210} follow from Corollary \ref{cor1} with the same argument as in \cite[Theorem 10.9]{GMR}, where the case $\alpha=\theta$ was considered. This gives the isometric behaviour of $\hat\Phi$. Moreover,
\begin{multline}\iota_2\circ\hat\Phi\circ\pi(\tfrac{\alpha}{\gcd(K_\theta, K_\alpha)}\varphi)= \hat\Phi(\tfrac{\alpha}{\gcd(K_\theta, K_\alpha)}\varphi+\alpha H^\infty)\\=A^{\theta,\alpha}_{\frac{\alpha}{\gcd(K_\theta, K_\alpha)}\varphi}=
\Phi(\tfrac{\alpha}{\gcd(K_\theta, K_\alpha)}\varphi)=
\Phi\circ\iota_1(\tfrac{\alpha}{\gcd(K_\theta, K_\alpha)}\varphi) . \end{multline}

Continuity of $\hat\Phi$ is a consequence of the fact that weak* topology on $\frac{\alpha}{\gcd(\alpha,\theta)}H^\infty\Big/ \alpha H^\infty$ is given by the quotient of weak* topology on  $\frac{\alpha}{\gcd(\alpha,\theta)}H^\infty$  as a restriction of weak* topology on $H^\infty$.

\end{proof}

In the case when $\alpha\leqslant\theta$ we have that
$$\hat{\Phi}\colon H^\infty\Big/ \alpha H^\infty\to \mathcal{I}(K_\theta, K_\alpha)$$
and $\Phi(\varphi)=\hat\Phi([\varphi])=\ata_\varphi$ for $\varphi\in H^\infty$. Hence by \cite[Theorem III.1.16]{Bercovici} and Theorem \ref{th72} we have the following:
\begin{corollary}\label{analityczne}
Let $\alpha$, $\theta$ be nonconstant inner functions such that $\alpha\leqslant\theta$. Then
\begin{enumerate}
\item $\mathcal{A}(K_\theta, K_\alpha)=\mathcal{I}(K_\theta, K_\alpha)$,%=\mathcal{W}(K_\theta, K_\alpha)
\item $\hat\Phi\colon H^\infty/\alpha H^\infty\to \mathcal{I}(K_\theta, K_\alpha)$ is a homeomorphism with respect to the corresponding weak* topologies,
\item The weak* topology and  the weak operator topology coincide on $\mathcal{I}(K_\theta, K_\alpha)$.
%\item $\mathcal{I}(K_\theta, K_\alpha)$ has the property $\mathbb{A}_1(1)$.
\end{enumerate}
 \end{corollary}

%\begin{proof}
%The proof is similar as in \cite[Proposition III.1.21]{Bercovici}.
%\end{proof}

%\begin{corollary}
%Let $\alpha$, $\theta$ be nonconstant inner functions such that $\alpha\leqslant\theta$.
%The space $\mathcal{A}(K_\theta, K_\alpha)$ is hereditarily 2-reflexive.
%\end{corollary}

%\begin{proof}
%Note firstly that Theorem \ref{analityczne}~(1) implies that $x\otimes S^*_\alpha y-S_\theta x\otimes y\in \mathcal{A}(K_\theta, K_\alpha)_\perp$ for all $x\in \kdt$, $y\in\kda$.
%
%On the other hand, if $A\in \mathcal{B}(\kdt,\kda)$ and $\operatorname{tr} (A(x\otimes S^*_\alpha y-S_\theta x\otimes y))=0$ for all $x\in \kdt$, $y\in\kda$, then
%$\langle Ax,S_\alpha^* y\rangle=\langle AS_\theta x, y\rangle.$ Hence $\langle S_\alpha A x, y\rangle=\langle AS_\theta x, y\rangle$, which implies that $A\in \mathcal{I}(K_\theta, K_\alpha)=\mathcal{A}(K_\theta, K_\alpha)$.
%\end{proof}

\section{Intertwining property for asymmetric dual truncated Toeplitz operators}
Let $\ta$ be nonconstant inner functions. For $\varphi \in L^2$ we can consider the densely defined multiplication operator
\begin{equation*}
  M_{\varphi}\colon (K_{\theta}\cap L^{\infty}) \oplus (K_{\theta}^{\bot}\cap L^{\infty}) \to K_{\alpha} \oplus K_{\alpha}^{\bot}.
  \end{equation*}
  According to this decomposition the operator $M_\varphi$
is given by the matrix $$\bmatrix P_{\alpha}M_{\varphi|K_{\theta}} & P_{\alpha}M_{\varphi|K_{\theta}^{\perp}} \\
P_{\alpha}^{\perp}M_{\varphi|K_{\theta}} & P_{\alpha}^{\perp}M_{\varphi|K_{\theta}^{\perp}}\endbmatrix =\bmatrix A^{\ta}_\varphi & (B^{\at}_{\bv})^* \\
B^{\ta}_{\varphi} & D^{\ta}_{\varphi}\endbmatrix .$$ %Here
%$P_\alpha {M_\varphi}_{|K_\theta^\perp}=(B^{\at}_{\bv})^*$.
%--------------------------------------------------------------------------------------
%
%$h\in K_{\theta}^{\bot}$, $g\in K_{\alpha}$
%\begin{equation*}
%  \langle (B^{\at}_{\bv})^*h,g\rangle = \langle h, B^{\at}_{\bv}g\rangle = \int h\varphi \bar g\, d\mu = \langle M_{\varphi}h,g\rangle
%\end{equation*}
%
%--------------------------------------------------------------------------------------

One can easily calculate that $M_\theta C_\alpha M_\varphi=M_\alpha M_{\bv}C_\theta$, that is \begin{equation}\label{caz}
C_\alpha M_\varphi=M_{\alpha \bv\bt}C_\theta.\end{equation}
 Writing that using matrices we get
\[\bmatrix {C_\alpha}_{|K_\alpha} & 0 \\
0 & {C_\alpha}_{|K_\alpha^\perp}\endbmatrix
\bmatrix A^{\ta}_\varphi & (B^{\at}_{\bv})^* \\
B^{\ta}_{\varphi} & D^{\ta}_{\varphi}\endbmatrix =
\bmatrix A^{\ta}_{\alpha\bv\bt} & (B^{\at}_{\ba\varphi\theta})^* \\
B^{\ta}_{\alpha\bv\bt} & D^{\ta}_{\alpha\bv\bt}\endbmatrix
\bmatrix {C_\theta}_{|K_\theta} & 0 \\
0 & {C_\theta}_{|K_\theta^\perp}\endbmatrix. \]

Thus we have the following
\begin{lemma}\label{L4.1} Let $\alpha$, $\theta$ be inner functions and $\varphi\in L^2$. Then
\begin{enumerate}\item $C_\alpha A^{\ta}_\varphi g=A^{\ta}_{\alpha\bv\bt} C_\theta g$ for $g\in K_\theta\cap L^\infty$;
\item $C_\alpha D^{\ta}_\varphi h=D^{\ta}_{\alpha\bv\bar\theta} C_\theta h$ for $h\in K_\theta^\perp\cap L^\infty$;
\item $C_\alpha (B^{\at}_{\bv})^*h=(B^{\at}_{\ba\varphi\theta})^*C_\theta h$ for $h\in K_\theta^\perp\cap L^\infty$;
\item $C_\alpha B^{\ta}_\varphi g=B^{\ta}_{\alpha\bv\bt}C_\theta g$ for $g\in K_\theta\cap L^\infty$.
\end{enumerate}
\end{lemma}
Recall that an asymmetric dual truncated Toeplitz operator is bounded if and only if it has a bounded symbol (see \cite[Proposition 1.1]{CKLPd}). We therefore have
$$\mathcal{T}(K^\perp_\theta, K^\perp_\alpha)=\{D_{\varphi}^{\theta,\alpha}\colon\, \varphi\in
L^\infty\}.$$
The following result describes $\mathcal{T}(K^\perp_\theta, K^\perp_\alpha)\cap \mathcal{I}(K^\perp_\theta, K^\perp_\alpha)$.

\begin{theorem}\label{idatto}
Let $\alpha,\theta$ be nonconstant inner functions and $\varphi \in L^{\infty}$, $\varphi\ne 0$. Then $D^{\ta}_{\varphi} D^{\theta}_{z}=D^{\alpha}_{z}D^{\ta}_{\varphi}$ if and only if one of the following holds
\begin{enumerate}
  \item $\alpha(0)=0=\theta(0)$ and $\varphi\in\frac{\alpha}{\mathrm{gcd}(\at)}K_{z\cdot \mathrm{gcd}(\at)}$, or
  \item $\alpha=\theta$ and $\varphi\in (\ktz)^{-1}K_{z\theta}$.
\end{enumerate}
\end{theorem}
\begin{proof}
Let $\varphi\in L^\infty\setminus\{0\}$. The commutation relation $M_\varphi M_z=M_z M_\varphi$ can be written as%we have on $L^2$,% a dense subset of $L^2$
\[\bmatrix A^{\ta}_{\varphi} & (B^{\at}_{\bv})^* \\
B^{\ta}_{\varphi} & D^{\ta}_{\varphi}\endbmatrix
\bmatrix A^{\theta}_z & (B^{\theta}_{\bar z})^* \\
B^{\theta}_z & D^{\theta}_z\endbmatrix =
\bmatrix A^{\alpha}_z & (B^{\alpha}_{\bar z})^* \\
B^{\alpha}_z & D^{\alpha}_z\endbmatrix
\bmatrix A^{\ta}_{\varphi} & (B^{\at}_{\bv})^*\\
B^{\ta}_{\varphi} & D^{\ta}_{\varphi}\endbmatrix. \]
Thus we have
\begin{equation*}B^{\ta}_{\varphi} (B^{\theta}_{\bar z})^*+D^{\ta}_{\varphi} D^{\theta}_{z}=B^{\alpha}_{z} (B^{\at}_{\bv})^*+D^{\alpha}_{z} D^{\ta}_{\varphi}.\end{equation*}
Clearly,%Note that by \eqref{caz}
 $$B^{\ta}_{\varphi}\ktz=P^{\bot}_{\alpha}(\varphi \ktz)$$
and, by Lemma \ref{L4.1}(4),
$$B^{\at}_{\bv}\tilde{k}^{\alpha}_0=B^{\at}_{\bv}C_{\alpha}\kaz=C_{\theta}B^{\at}_{\theta \varphi \ba}\kaz=C_{\theta}P^{\bot}_{\theta}(\theta \varphi\ba \kaz).$$ Therefore, by \cite[Corollary 7]{CKLPcom}, we get
\begin{equation*}
\aligned
D^{\ta}_{\varphi}D^{\theta}_{z}-D^{\alpha}_{z} D^{\ta}_{\alpha}&=B^{\alpha}_{z}(B^{\ta}_{\varphi})^*-B^{\ta}_{\varphi}(B^{\theta}_{\bar z})^*\\
&=(\alpha \otimes \tilde{k}^{\alpha}_0)(B^{\at}_{\bv})^*-B^{\ta}_{\varphi}(k^{\theta}_0 \otimes \bar z) \\
&=\alpha \otimes (B^{\at}_{\bv} \tilde{k}^{\alpha}_0)-(B^{\ta}_{\varphi}k^{\theta}_0)\otimes \bar z\\
%&=\alpha \otimes C_{\theta}P^{\bot}_{\theta}(\theta\ba \varphi k^{\alpha}_0)-P^{\bot}_{\alpha}(\varphi k^{\theta}_0)\otimes \bar z\\
&=\alpha \otimes C_{\theta}P^{\bot}_{\theta}(\theta\ba \varphi k^{\alpha}_0)-P^{\bot}_{\alpha}(\varphi k^{\theta}_0)\otimes C_{\theta}\theta
\endaligned
\end{equation*}
(compare with \cite[Theorem 11(1)]{CKLPcom}). Hence
\begin{equation}
D^{\ta}_{\varphi} D^{\theta}_{z}=D^{\alpha}_{z} D^{\ta}_{\varphi}\end{equation}
if and only if there is a constant $c\in \mathbb{C}$ such that
\begin{equation}
\left\{ \aligned P^{\bot}_{\alpha}(\varphi k^{\theta}_0)&=c{\,\alpha}, \\ P^{\bot}_{\theta}(\theta \ba \varphi k^{\alpha}_0)&=c\,{\theta}.\endaligned\right.\end{equation}
Equivalently, there are $g\in K_\alpha$ and $h\in K_\theta$ such that
\begin{align}\label{cg}\varphi \ktz&=c\alpha+g \qquad \text{and}\\ %\Longleftarrow \varphi \in H^2$
\label{ch}\theta\ba\varphi\kaz&=c\theta+h.\end{align}
Since the functions $k_0^\alpha$, $k_0^\theta$ are bounded from below and  analytic (i.e., $(k_0^\alpha)^{-1}, (k_0^\theta)^{-1}\in H^\infty$), we get by \eqref{cg} that $\varphi \in H^2$.
%$B^{\ta}_{\varphi}(B^{\theta}_{\bar z})^*+D^{\ta}_{\varphi}D_z=B^{\alpha}_{z}(B^{\at}_{\bv})^*+D^{\alpha}_{z}D^{\ta}_{\varphi}$
%
%$D^{\ta}_{\varphi}D_z-D^{\alpha}_{z}D^{\ta}_{\varphi}=(\alpha \otimes \tilde{k}^{\alpha}_0)(B^{\at}_{\bv})^*-B^{\ta}_{\varphi}(\ktz\theta \bar z)=\alpha \otimes (B^{\at}_{\bv}\tilde{k}^{\alpha}_0) - (B^{\ta}_{\varphi}\ktz) \otimes \bar z = \alpha \otimes C_{\theta}P^{\bot}_{\theta}(\theta\ba \varphi \kaz)P^{\bot}_{\alpha}(\varphi \ktz)\otimes C_{\theta}\theta$
Let $\gamma=\text{gcd}(\alpha,\theta)$. Then, by \eqref{ch},
$$\varphi\kaz-c\alpha=\bt\alpha h= \tfrac\bt\bg \tfrac\alpha\gamma h \in H^2.$$
Hence $h$ is divisible by $\frac{\theta}{\gamma}$ and
since $h\in K_{\theta}=K_{\frac\theta\gamma}\oplus \frac\theta\gamma K_{\gamma}$, we have
%$\frac\alpha\gamma h \in \frac\theta\gamma H^2$
%
$$h=\tfrac\theta\gamma h_1\quad \text{with }h_1\in K_{\gamma}.$$
Therefore, by \eqref{ch},
$$\tfrac\ba\bg \varphi \kaz=c \gamma +  h_1\in H^2.$$
Since $\kaz$ is an outer function, it cannot be divisible by $\frac\alpha\gamma$, which implies that
$$\varphi=\tfrac\alpha\gamma \varphi_1 \quad\text{with }\varphi_1\in H^2\setminus\{0\}.$$
Hence by \eqref{cg} we get
$\frac\alpha\gamma\varphi_1\ktz=c\alpha+g$. Therefore  $g=\frac\alpha\gamma g_1$, $g_1\in K_{\gamma}$, and \eqref{cg}, \eqref{ch} are equivalent to
\begin{align}\label{cg1}\varphi_1\ktz&=c\gamma+g_1, \quad g_1\in K_{\gamma},\\ \label{ch1}
\varphi_1\kaz&=c\gamma+h_1,\quad h_1\in K_{\gamma}.\end{align}
Moreover,
$$P_\gamma(\varphi_1)=P_\gamma(\varphi_1\ktz+\overline{\theta(0)}\theta\varphi_1)=
P_\gamma(c\gamma+g_1+\overline{\theta(0)}\theta\varphi_1)=g_1.$$
Similarly,
$P_\gamma(\varphi_1)=h_1$, so $g_1=h_1$. Comparing \eqref{cg1} with \eqref{ch1} we thus get
$\varphi_1(\kaz-\ktz)=0$. Since $\varphi_1,(\kaz-\ktz)\in H^2$ and
$\varphi_1\neq0$, we must have $$\kaz-\ktz=\overline{\theta(0)}\theta-\overline{\alpha(0)}\alpha=0.$$
This is possible only in two cases:
\begin{enumerate}

  \item if $\alpha(0)=0=\theta(0)$, then $k_0^\alpha=k_0^\theta=1$ and by \eqref{cg1},
  $$\varphi=\tfrac{\alpha}{\gamma}\varphi_1\in \tfrac{\alpha}{\gamma}(K_\gamma\oplus \mathbb{C}\gamma)= \tfrac{\alpha}{\gamma}K_{z \cdot \gamma};$$
   \item if $\alpha=\theta$, then \eqref{cg1} and \eqref{ch1} become the same condition, %are carbon copy and
   equivalent to $\varphi\ktz\in K_\theta\oplus \mathbb{C}\theta=K_{z\theta}$, which leads to $\varphi\in (\ktz)^{-1}K_{z\theta}$.
\end{enumerate}
\end{proof}

The following Corollary was recently obtained independently in \cite{Gu2} (see \cite[Corollary 2.6]{Gu2}).

\begin{corollary}\label{wnios}
	Let $\theta$ be a nonconstant inner function and let $\varphi \in L^{\infty}$, $\varphi\ne 0$. Then $D^{\theta}_{\varphi}\in\{D^{\theta}_{z}\}'$ if and only if $\varphi\in (\ktz)^{-1}K_{z\theta}$.
\end{corollary}

\section{Characterization of operators intertwining $D^\theta_z$ and $D^\alpha_z$}

In this section we will characterize operators from $\mathcal{I}(K_\theta^\perp, K_\alpha^\perp)$. We will use the description of shift invariant operators from $\mathcal{B}(K_\theta^\perp, K_\alpha^\perp)$, given in \cite{OurNewPaper}.

Recall that $D\in \mathcal{B}(K_\theta^\perp, K_\alpha^\perp)$ is called {\it shift invariant} if \begin{equation*}%\label{shin}
\langle D M_z f,M_z g\rangle=\langle D f,g\rangle%\quad
%\text{for}\  f\in K_\theta^\perp\cap\{\bar z\}^\perp,\,  g\in K_\alpha^\perp\cap\{\bar z\}^\perp.
\end{equation*}
for all $f\in K_\theta^\perp$ and $g\in K_\alpha^\perp$ such that $M_zf\in K_\theta^\perp$ and $M_zg\in K_\alpha^\perp$, that is, $f\in K_\theta^\perp\cap\{\bar z\}^\perp$ and $g\in K_\alpha^\perp\cap\{\bar z\}^\perp$ (see \cite{OurNewPaper}). The characterization of shift invariant elements of $\mathcal{B}(K_\theta^\perp, K_\alpha^\perp)$ is expressed in terms of compressions of $M_\varphi$ to some subspaces of $K_\theta^\perp$. Using decompositions  $K_\theta^\perp=\theta H^2\oplus H^2_- $ and $K_\alpha^\perp=\alpha H^2\oplus H^2_- $ one can write each operator $D\in \mathcal{B}(K_\theta^\perp, K_\alpha^\perp)$ as a matrix
$$D=\begin{bmatrix}\dtha&\dba\\ \dc &\dd\end{bmatrix}.$$
In particular, for $\varphi\in L^{\infty}$, we obtain
$$D_\varphi^{\theta,\alpha}=\begin{bmatrix}\hat{T}_\varphi^{\theta,\alpha}&\check{\Gamma}_{\varphi}^\alpha\\ \hat{\Gamma}_\varphi^\theta &\check{T}_\varphi\end{bmatrix}=\begin{bmatrix}\hat{T}_\varphi^{\theta,\alpha}&(\hat{\Gamma}_{\bar{\varphi}}^{{\alpha}})^*\\
\hat{\Gamma}_\varphi^{\theta}&\check{T}_\varphi\end{bmatrix},$$
where
\begin{align*} \hat{T}_\varphi^{\theta,\alpha}&=P_{\alpha H^2}M_{\varphi|\theta H^2},  \quad \check{\Gamma}_\varphi^{\alpha}=P_{\alpha H^2}M_{\varphi| H^2_-}, \\  \hat{\Gamma}_\varphi^{\theta}&=P^-M_{\varphi|\theta H^2}\quad\text{and}\quad  \check{T}_\varphi=P^-M_{\varphi| H^2_-}.
\end{align*}
Clearly, each of the operators above can be defined for $\varphi\in L^2$ (densely, for bounded functions).

We will now consider operators of the form
	\begin{equation*}
%\label{tmat}
D=\begin{bmatrix}\hat{T}_{\varphi_1}^{\theta,\alpha}&\check{\Gamma}_{\varphi_2}^{\alpha}
\\\hat{\Gamma}_{\varphi_3}^{\theta}&\check{T}_{\varphi_4}\end{bmatrix}.
\end{equation*}
with $\varphi_i\in L^2$ for $i=1,2,3,4$. Note that if the operator $D$ given above is bounded, then necessarily $\varphi_1,\varphi_4\in L^\infty$. This is a consequence of the fact that $\hat{T}_{\varphi_1}^{\theta,\alpha}$ and $\check{T}_{\varphi_4}$ are determined by the classical Toeplitz operators $T_{\varphi_1}$ and $T_{\bar \varphi_4}$, respectively (see \cite[Proposition 20]{CKLPcom}), and that a classical Toeplitz operator is bounded if and only if its symbol is from $L^{\infty}$. On the other hand, even if $D$ is bounded, the functions $\varphi_2$ and $\varphi_3$ may not belong to $L^\infty$ (since a classical Hankel operator may be bounded even if its symbol is not, see \cite[Chapter 1]{VVP} for details).

\begin{lemma}\label{l6.1} Let $\alpha,\theta $ be nonconstant inner functions. Then
\begin{enumerate}[{(i)}]
\item $\hat{T}_{\varphi_1}^{\theta,\alpha}=0$ if and only if $\varphi_1=0$;
\item $\check{\Gamma}_{\varphi_2}^{\alpha}=0$ if and only if $\varphi_2\perp \alpha z H^2$, that is, $\varphi_2\in\alpha \overline{H^2}$. In particular, $\check{\Gamma}_{\varphi_2}^{\alpha}=0$ for $\varphi_2\in\alpha \overline{H^\infty}$;%$\varphi_2\in\alpha \overline{H^\infty}$.
\item $\hat{\Gamma}_{\varphi_3}^{\theta}=0 $ if and only if $\varphi_3\perp \overline{\theta z H^2}$, that is, $\varphi_3\in \bar \theta H^2$. In particular, $\hat{\Gamma}_{\varphi_3}^{\theta}=0 $ for $\varphi_3\in \bar \theta H^\infty$;%$\varphi_3\in \bar \theta H^\infty$.
\item $\check{T}_{\varphi_4}=0$ if and only if $\varphi=0$.
\end{enumerate}
\end{lemma}
\begin{proof}
See \cite{CKLPcom}.
\end{proof}

%\begin{lemma}\label{l8.1}
%	Let $f\in K_\theta^\perp$. Then $M_z f\in K_\theta^\perp$ if and only if $f$ is orthogonal to the function $\bar z\in L^2$.
%\end{lemma}

%\begin{proof}
%	Let $f=\bar z \bar g+\theta h\in K_\theta^\perp$, $g, h\in H^2$. Note that $M_z f=\bar g+z\theta h\in K_\theta^\perp$ only if
%	$$0=\langle\bar g, k_0^\theta \rangle=\overline{g(0)}(1-|\theta (0)|^2).$$ Hence $\overline{g(0)}=0$ and
%	$g=zg_1$, for some $g_1\in H^2$.
%\end{proof}

\begin{theorem}[\cite{OurNewPaper}]\label{shinth}
	Let $\theta,\alpha$ be two nonconstant inner functions and let $K_\theta^\perp=\theta H^2 \oplus H^2_-$, $K_\alpha^\perp=\alpha H^2 \oplus H^2_-$. The operator $D\in \mathcal{B}(K_\theta^\perp,K_\alpha^\perp)$ is shift invariant if and only if there are $\varphi_i\in L^\infty$ for $i=1,2,3,4$ such that
	\begin{equation*}
	%\label{tmat}
	D=\begin{bmatrix}\hat{T}_{\varphi_1}^{\theta,\alpha}&\check{\Gamma}_{\varphi_2}^{\alpha}
	\\\hat{\Gamma}_{\varphi_3}^{\theta}&\check{T}_{\varphi_4}\end{bmatrix}.
	\end{equation*}
\end{theorem}

We can now state our main result.

\begin{theorem}\label{kmutant}
Let $\alpha,\theta $ be nonconstant inner functions. The operator $D\in \mathcal{B}(K_\theta^\perp,K_\alpha^\perp)$ intertwines $D_z^{\theta}$ and $D_z^{\alpha}$, i.e.,
\begin{equation}\label{inttwin} DD_z^{\theta}= D_z^{\alpha}D
\end{equation}
 if and only if $D\in
\mathcal{T}^2(K_\theta^\perp,K_\alpha^\perp)$, i.e., $D=\begin{bmatrix}\hat{T}_{\psi_1}^{\theta,\alpha}&\check{\Gamma}_{\psi_2}^{\alpha}
\\\hat{\Gamma}_{\psi_3}^{\theta}&\check{T}_{\psi_4}\end{bmatrix}$, where the functions $\psi_i\in L^2$, $ i=1,2,3,4$, satisfy the following conditions:
\begin{enumerate}[{(i)}]
\item if  $\theta(0)\ne 0$, then
    \begin{equation}\label{c1c2n}\begin{split}
\psi_1&=\tfrac {\overline{\alpha(0)}} {\overline{\theta(0)}}\alpha\bar\theta \psi_4,
  \qquad \psi_2= {{\overline{\alpha(0)}}}\alpha P^+( \psi_4),\\
 \psi_3&= \tfrac{1}{{\overline{\theta(0)}}}\,\bar\theta\, P^-( \psi_4),\qquad \psi_4\in L^\infty .
 \end{split}\end{equation}
%
%\item if $\alpha(0)= 0$, $\theta(0)\ne 0$,  then  $\psi_4\in L^\infty$ are arbitrary and
%   \begin{equation}\label{fi36}
%\psi_1=0,\quad\
%  \psi_2= 0,\quad\
% \psi_3= \tfrac{1}{{\overline{\theta(0)}}}\,\bar\theta\, P^-( \psi_4). \end{equation}
\item if ${{\theta(0)}}=0$ and $\alpha(0)\ne 0$, then
   \begin{equation}\label{fi32}\psi_1\in L^\infty,\quad
  \psi_2= 0,\quad
 \psi_3= \tfrac{1}{{\overline{\alpha(0)}}}\,\bar\theta\, P^-(\bar \alpha \theta \psi_1),\quad\psi_4=0. \end{equation}
     \item if ${{\theta(0)}}={{\alpha(0)}}=0$, then  \begin{equation}\label{fi33}\psi_1\in \alpha\bar\theta H^\infty,\quad\psi_2=0,\quad \psi_3\in L^{\infty},\quad \psi_4 \in H^\infty.\end{equation}
\end{enumerate}
\end{theorem}

Corollary \ref{5.3} should be compared with \cite[Proposition 3.1]{Gu2} and \cite[Theorem 1.1]{LSD}.
\begin{corollary}\label{5.3} Let $\theta $ be a nonconstant inner function. The operator $D\in \mathcal{B}(K_\theta^\perp)$ commutes with  $D_z^{\theta}$
 if and only if $D\in
\mathcal{T}^2(K_\theta^\perp)$, i.e., $D=\begin{bmatrix}\hat{T}_{\psi_1}^{\theta}&\check{\Gamma}_{\psi_2}^{\theta}
\\ \hat{\Gamma}_{\psi_3}^{\theta}&\check{T}_{\psi_4}\end{bmatrix}$, where the functions $\psi_i\in L^2$, $ i=1,2,3,4$, satisfy the following conditions:
\begin{enumerate}[{(i)}]
\item if $\theta(0)\ne 0$, then
       \begin{equation}\label{c1c2na}
\psi_1= \psi_4,
  \ \psi_2= {{\overline{\theta(0)}}}\theta P^+(\psi_4),
 \ \psi_3= \tfrac{1}{{\overline{\theta(0)}}}\,\bar\theta\, P^-( \psi_4),\ \psi_4\in L^\infty
  . \end{equation}
%\end{equation*}
\item if $ \theta(0)= 0$, then \begin{equation}\psi_1\in  H^\infty,\quad\psi_2=0,\quad\psi_3\in L^\infty,\quad \psi_4 \in H^\infty.\end{equation}
\end{enumerate}
\end{corollary}

\begin{proof}[Proof of Theorem \ref{kmutant}] Let $D\in \mathcal{B}(K_\theta^\perp,K_\alpha^\perp)$. We first show that if $D$ satisfies \eqref{inttwin}, then it is shift invariant. To this end, take $f\in K_\theta^\perp$ and $g\in K_\alpha^\perp$ such that $M_zf\in K_\theta^\perp$ and $M_zg\in K_\alpha^\perp$.  Then $D_z^{\theta}f=M_z f$, $D_z^{\alpha}g=M_zg$ and so
\begin{align*}%\label{shin2}
  \langle D M_z f,M_z g\rangle&=\langle D D^{\theta}_z f,M_z g\rangle=\langle D_z^{\alpha}D  f,M_z g\rangle\\&=\langle  M_z D f,M_z g\rangle=\langle D f,g\rangle,
\end{align*}  i.e., $D$ is shift invariant. Therefore, by Theorem \ref{shinth}, to describe all $D$ satisfying \eqref{inttwin} it is enough to consider operators of the form
 \begin{equation}
\label{tmat}D=\begin{bmatrix}\hat{T}_{\varphi_1}^{\theta,\alpha}&\check{\Gamma}_{\varphi_2}^{\alpha}
\\\hat{\Gamma}_{\varphi_3}^{\theta}&\check{T}_{\varphi_4}\end{bmatrix},
\end{equation}
where $\varphi_i\in L^\infty$ for $i=1,2,3,4$. Assume then that $D$ is given by \eqref{tmat}.
% Consider $D_z^{\alpha}$ then
Recall that $\check{\Gamma}_z^{\alpha} \colon H^2_- \to \alpha H^2$ and for $h\in z H^2$ we have
\begin{align*}
\check{\Gamma}^\alpha_z \bar h&= P_{\alpha H^2}(z\bar h)=P_{\alpha H^2}P^+(z\bar h)\\&=P_{\alpha H^2}(\langle \bar h, \bar z\rangle \cdot 1)= \overline {\alpha(0)}\langle \bar h, \bar z\rangle \alpha=\overline {\alpha(0)}\ (\alpha\otimes \bar z) \bar h.
\end{align*}
On the other hand,
$\hat{\Gamma}_z^{\alpha} \colon \alpha H^2 \to H^2_-$ and for $h\in H^2$ we have
\[\hat{\Gamma}_z^{\alpha}\alpha h=P^-(z \alpha h)=0.
 \] Similarly, $\check{\Gamma}^\theta_z=\overline {\theta(0)}\ (\theta\otimes \bar z)$ and $\hat{\Gamma}_z^{\theta}=0$. % calculations we make for  $\theta$. Therefore
%${\overline{\theta(0)}}=\overline {\theta(0)}$  $c_2=\overline {\alpha(0)}$
It follows that
\begin{equation}\label{dz}
D^{\theta}_z=\bmatrix \hat{T}^\theta_z & {\overline{\theta(0)}}(\theta \otimes \bar z) \\
0 & \check{T}_z\endbmatrix \quad \text{and}\quad D^{\alpha}_z=
\bmatrix \hat{T}^{\alpha}_z & {\overline{\alpha(0)}}(\alpha\otimes\bar z) \\
0 & \check{T}_z \endbmatrix.
\end{equation}
Therefore the intertwining relation \eqref{inttwin} can be expressed as
\begin{equation*}%\label{commut}
\bmatrix \hat{T}^{\ta}_{\varphi_1} & \check{\Gamma}^{\alpha}_{\varphi_2} \\
\hat{\Gamma}^{\theta}_{\varphi_3} & \check{T}_{\varphi_4} \endbmatrix
\bmatrix \hat{T}^\theta_z & {\overline{\theta(0)}}(\theta \otimes \bar z) \\
0 & \check{T}_z\endbmatrix =
\bmatrix \hat{T}^{\alpha}_z & {\overline{\alpha(0)}}(\alpha\otimes\bar z) \\
0 & \check{T}_z \endbmatrix
\bmatrix \hat{T}^{\ta}_{\varphi_1} & \check{\Gamma}^{\alpha}_{\varphi_2}\\
\hat{\Gamma}^{\theta}_{\varphi_3} & \check{T}_{\varphi_4}\endbmatrix,
%\bmatrix \theta h \\ \bar g \endbmatrix
\end{equation*}
or equivalently as the following set of conditions%Comparing both sides we get
\begin{align}
&\label{11} &\hat{T}^{\ta}_{\varphi_1}\hat{T}^{\theta}_z  &=\hat{T}^{\alpha}_z\hat{T}^{\ta}_{\varphi_1}+{\overline{\alpha(0)}}(\alpha\otimes\bar z) \ \hat\Gamma^{\theta}_{\varphi_3};
\\
&\label{12} &{\overline{\theta(0)}}\ \hat{T}^{\ta}_{\varphi_1}(\theta)\otimes\bar z + \check\Gamma^{\alpha}_{\varphi_2}\check{T}_z&=\hat{T}^{\alpha}_z\check \Gamma^{\alpha}_{\varphi_2}+{\overline{\alpha(0)}}\ (\alpha\otimes\bar z)\check{T}_{\varphi_4};
\\
&\label{21}   &\hat\Gamma^{\theta}_{\varphi_3}\hat{T}^{\theta}_z&=\check{T}_z\hat\Gamma^{\theta}_{\varphi_3};
\\
&\label{22} &{\overline{\theta(0)}}\ \hat\Gamma^{\theta}_{\varphi_3}          (\theta)\otimes \bar z + \check{T}_{\varphi_4}\check{T}_z&=\check{T}_z\check{T}_{\varphi_4}.
\end{align}
%
%$P^-_{\varphi_3}$ $z\theta h=P^-zP^-\varphi_3\theta h=0$
%
%-----------------------------------------
%
%(1) $\theta,\alpha$ $x\in \mathcal{I}(\ta)\Longrightarrow \exists_{u\in H^{\infty}}: \alpha|_{u\theta},\alpha\in H^2$, $\|u\|=\|x\|$, $a\in A_u$
%
%(2) $\forall_{u\in H^{\infty}}: \alpha|_{u\theta} \,x=A^{\ta}_u\in \mathcal{I}(\ta)$ $x=0 \Longleftrightarrow \alpha_{u}$
%
%
%\newpage
Note firstly that  the condition \eqref{22} is equivalent to
\begin{equation}
\label{plus}\check{T}_z\check T_{\varphi_4}-\check{T}_{\varphi_4}\check{T}_z={\overline{\theta(0)}}\ \hat\Gamma^{\theta}_{\varphi_3}(\theta)\otimes \bar z.
\end{equation}
For $ h\in zH^2$ we have
\begin{align*}
(\check{T}_z\check T_{\varphi_4}-\check{T}_{\varphi_4}\check{T}_z)\bar h&=P^-(zP^-(\varphi_4\bar h))-P^-(\varphi_4P^-(z\bar h))\\
&=P^-(zP^-(\varphi_4\bar h))-P^-(\varphi_4(z\bar h- \langle \bar h, \bar z \rangle\cdot 1))\\
&=P^-(zP^-(\varphi_4\bar h))-P^-(\varphi_4z\bar h)+\langle \bar h,\bar z\rangle P^-(\varphi_4)\\
&=P^-(z(P^-(\varphi_4\bar h)-\varphi_4\bar h))+\langle \bar h,\bar z\rangle P^-(\varphi_4).
\end{align*}
%hence \eqref{22} is equivalent to
%$$P^-(zP^-(\varphi_4\bar h))-P^-\varphi_4z\bar h+\langle \bar h,\bar z\rangle P^-(\varphi_4\bar z )=-{\overline{\theta(0)}} \langle \bar h, \bar z\rangle  P^-(\varphi_3\theta).$$
%
%$=-{\overline{\theta(0)}} \langle \bar h,\bar z\rangle P^-(\varphi_3\theta)$
%
%$=-{\overline{\theta(0)}} \langle \bar h,\bar z\rangle P^-(\varphi_3\theta)$
Because $P^--I_{L^2}=-P^+$, thus
$$P^-(z(P^-(\varphi_4\bar h)-\varphi_4\bar h))=-P^-(zP^+(\varphi_4\bar h))=0$$
and \eqref{plus} can be written as
\begin{equation*}
P^-(\varphi_4)\otimes \bar z=P^-(\overline{\theta(0)}\theta\varphi_3)\otimes \bar z,
\end{equation*}
%$$+\langle \bar h,\bar z\rangle P^-(\varphi_4)={\overline{\theta(0)}} \langle \bar h,\bar z\rangle P^-(\theta\varphi_3).$$
%$$-P^-(z(P^+(\varphi_4\bar h))={\overline{\theta(0)}} \langle \bar h,\bar z\rangle P^-(\varphi_3\theta)-\langle \bar h,\bar z\rangle P^-(\varphi_4).$$
or
$$P^-(\varphi_4-{\overline{\theta(0)}}\theta\varphi_3)\otimes\bar z=0,$$
and finally \eqref{22} is equivalent to
\begin{equation} P^-(\varphi_4-{\overline{\theta(0)}}\theta\varphi_3)=0.\end{equation}

Note that \eqref{21} is always satisfied. Indeed, for $h\in H^2$,
%(2.1) $\hat\Gamma^{\theta}_{\varphi_3}\hat{T}^{\theta}_{z}=\check{T}_z\hat\Gamma^{\theta}_{\varphi_3}$ always true
\begin{align*}
\hat\Gamma^{\theta}_{\varphi_3}\hat{T}^{\theta}_z(\theta h)&=
P^-({\varphi_3}P_{\theta H^2}(z\theta h))=
P^-({\varphi_3}z\theta h)\\
&=P^-(zP^-({\varphi_3}\theta h)+zP^+({\varphi_3}\theta h))\\
&=P^-(zP^-({\varphi_3}\theta h))=\check{T}_z\hat\Gamma^{\theta}_{\varphi_3}(\theta h).
\end{align*}
%if and only if
%$$P^-(\varphi_3z\theta h-zP^-(\varphi_3\theta h))=0.$$ That is equivalent to
%$$P^-(z(\varphi_3\theta h-P^-(\varphi_3\theta h))=0,$$ i.e.,
%$P^-(zP^+(\varphi_3\theta h))=0$, which is always true.

Let us now focus on condition \eqref{11}. Note that for $ h\in H^2$,
%$$\hat{T}^{\ta}_{\varphi_1}\hat T^{\theta}_z=
%\hat{T}^{\alpha}_z\hat{T}^{\ta}_{\varphi_1}-{\overline{\alpha(0)}}(\alpha\otimes\bar z) \hat\Gamma^{\theta}_{\varphi_3}$$
%and
\begin{equation}\begin{split}\label{g11}(\alpha\otimes\bar z) \hat\Gamma^{\theta}_{\varphi_3}(\theta h)&=(\alpha\otimes (\hat\Gamma^{\theta}_{\varphi_3})^*{\bar z})(\theta h)\\&=(\alpha\otimes (\check\Gamma^{\theta}_{\bar\varphi_3}
\bar z)) (\theta h)=\langle \theta h,\bar z\bar\varphi_3\rangle\alpha,\end{split}\end{equation}
and
\begin{multline*}(\hat{T}^{\ta}_{\varphi_1}\hat T^{\theta}_z-
\hat{T}^{\alpha}_z\hat{T}^{\ta}_{\varphi_1})(\theta h)=P_{\alpha H^2}(\varphi_1P_{\theta H^2}(z\theta h))-P_{\alpha H^2}(zP_{\alpha H^2}(\varphi_1\theta h))\\
=P_{\alpha H^2}(z(\varphi_1\theta h-P_{\alpha H^2}(\varphi_1\theta h)))=P_{\alpha H^2}(z(P^-+P_{\alpha})(\varphi_1\theta h)).
\end{multline*}
Since
$P^-+P_\alpha= \alpha P^-\bar\alpha$, \eqref{g11} and the calculations above show that \eqref{11} is equivalent to\begin{equation}
\label{plusplus}P_{\alpha H^2}(z\alpha P^-(\bar\alpha\theta\varphi_1h))={\overline{\alpha(0)}}\langle z\theta h,\bar\varphi_3\rangle\alpha\quad\text{for all }h\in H^2.
\end{equation}
%Hence
%$$\langle \bar\alpha\theta\varphi_1 h,\bar z\rangle\alpha ={\overline{\alpha(0)}}\langle\theta h,\bar z\bar\varphi_3\rangle\alpha$$
%and so
Since
\begin{align*}
P_{\alpha H^2}(z\alpha P^-(\bar\alpha\theta\varphi_1h))&=%\alpha
%P^+(\bar{\alpha}z\alpha P^-(\bar\alpha\theta\varphi_1h))\\
%&=
\alpha
P^+(z P^-(\bar\alpha\theta\varphi_1h))\\&=\alpha \langle zP^-(\bar\alpha\theta\varphi_1 h),1\rangle=\langle \bar\alpha\theta\varphi_1 h,\bar z\rangle\alpha,
\end{align*}
condition \eqref{plusplus} can be expressed as
$$\langle\bar\alpha\theta\varphi_1 -{\overline{\alpha(0)}}\theta\varphi_3,\bar z\bar h\rangle=0.$$
We have thus proved that \eqref{11} is equivalent to $$P^-(\bar\alpha\theta\varphi_1-{\overline{\alpha(0)}}\theta\varphi_3)=0\quad\text{for all }h\in H^2.$$
%
%-----------------------------------------------------------------------------
%
%(1.2 )
%$-{\overline{\theta(0)}}\ \hat{T}^{\ta}_{\varphi_1}(\theta)\otimes\bar z + \check\Gamma^{\alpha}_{\varphi_2}\check{T}_z=\hat{T}^{\alpha}_z\check \Gamma^{\alpha}_{\varphi_2}-{\overline{\alpha(0)}}\ (\alpha\otimes\bar z)\check{T}_{\varphi_4}$

Finally, condition \eqref{12} is equivalent to
$$\check\Gamma^{\alpha}_{\varphi_2}\check{T}_z - \hat{T}^{\alpha}_z\check\Gamma^{\alpha}_{\varphi_2}={\overline{\alpha(0)}}(\alpha \otimes \check{T}_{\varphi_4}^*(\bar z))-{\overline{\theta(0)}}(\hat{T}^{\ta}_{\varphi_1}(\theta)\otimes \bar z).$$
Now for any $h\in zH^2$ we have
\begin{align*}
(\check\Gamma^{\alpha}_{\varphi_2}\check{T}_z - \hat{T}^{\alpha}_z\check\Gamma^{\alpha}_{\varphi_2})\bar h&=P_{\alpha H^2}(\varphi_2P^-(z \bar h))-P_{\alpha H^2}(zP_{\alpha H^2}(\varphi_2\bar h))\\
%&={\overline{\alpha(0)}}\langle \bar h, P^-(\bar\varphi_4\bar z)\rangle \alpha-{\overline{\theta(0)}} \langle \bar h,\bar z\rangle P_{\alpha H^2}(\theta\varphi_1).
&=P_{\alpha H^2}(\varphi_2(z \bar h - \langle \bar h,\bar z\rangle\cdot 1))-zP_{\alpha H^2}(\varphi_2\bar h)\\
&=P_{\alpha H^2}(z\varphi_2 \bar h)-zP_{\alpha H^2}(\varphi_2\bar h)-\langle \bar h,\bar z\rangle P_{\alpha H^2}(\varphi_2).
\end{align*}
%Hence
%\begin{multline}\label{l12}P_{\alpha H^2}(\varphi_2(z \bar h - \langle \bar h,\bar z\rangle\cdot 1))-zP_{\alpha H^2}(\varphi_2\bar h)\\=({\overline{\alpha(0)}} (\alpha \otimes P^-(\bar z\bar\varphi_4))-{\overline{\theta(0)}}(P_{\alpha H^2}(\theta\varphi_1)\otimes \bar z))\bar h.\end{multline}
Since for each $f\in L^2$,
 \begin{align*}P_{\alpha H^2}(zf)- zP_{\alpha H^2}f&=\alpha P_+(z\bar\alpha f)-\alpha zP_+(\bar\alpha f)\\
 &= \alpha\langle \bar\alpha f,\bar z\rangle\cdot 1=\langle f,\alpha \bar z\rangle\alpha,\end{align*}
 it follows that
%thus the left hand side of \eqref{l12}  is equal to
 \begin{align*}
 (\check\Gamma^{\alpha}_{\varphi_2}\check{T}_z - \hat{T}^{\alpha}_z\check\Gamma^{\alpha}_{\varphi_2})\bar h&=\langle \varphi_2\bar h, \bar z\alpha\rangle\alpha-\langle \bar h,\bar z\rangle P_{\alpha H^2}(\varphi_2)\\
 &=(\alpha \otimes P^-(\bar z\alpha\bar\varphi_2)-P_{\alpha H^2}(\varphi_2)\otimes \bar z)\bar h.\end{align*}
Thus \eqref{12} is equivalent to
$$\alpha \otimes P^-(\bar z\alpha\bar\varphi_2)-P_{\alpha H^2}(\varphi_2)\otimes \bar z={\overline{\alpha(0)}}(\alpha \otimes \check{T}_{\varphi_4}^*(\bar z))-{\overline{\theta(0)}}(\hat{T}^{\ta}_{\varphi_1}(\theta)\otimes \bar z),$$
which in turn is equivalent to
\begin{equation}
\label{dwao}
 \alpha \otimes P^-(\bar z\alpha\bar\varphi_2- {{\alpha(0)}}\bar z\bar\varphi_4)=P_{\alpha H^2}(\varphi_2-{\overline{\theta(0)}}\theta\varphi_1)\otimes \bar z.
 \end{equation}
 Now \eqref{dwao} is possible only if for some $c\in \mathbb{C}$,
$$\left\{\aligned P_{\alpha H^2}(\varphi_2-{\overline{\theta(0)}}\theta\varphi_1)=c\alpha,\\
P^-(\bar z\alpha\bar\varphi_2- {{\alpha(0)}}\bar z\bar\varphi_4)=\bar c\bar z,\endaligned\right.$$
or equivalently
\begin{equation}\label{aa}
\left\{\aligned P^+(\bar{\alpha}\varphi_2-{\overline{\theta(0)}}\bar{\alpha}\theta\varphi_1)=c,\\
P^+(\bar \alpha\varphi_2- \overline{{\alpha(0)}}\varphi_4)=c.\endaligned\right.
\end{equation}
The second condition follows from the fact that $P^-(\bar z f)=\bar z \overline{P^+(\bar f)}$ for each $f\in L^2$.

%The second condition is equivalent to
%$P^+(\varphi_2\bar\alpha- {\overline{\alpha(0)}}\varphi_4)= c$.
Summing up, we have shown that the intertwining relation \eqref{inttwin} is equivalent to the following set of conditions%
%\begin{equation}\label{c12}
%\left\{
%\aligned
%P^-(\bar\alpha\theta\varphi_1-{\overline{\alpha(0)}}\theta\varphi_3)&=0,\\
%P^-(\varphi_4-{\overline{\theta(0)}}\theta\varphi_3)&=0,
%\\
% P^+(\bar\alpha\varphi_2-{\overline{\theta(0)}}\bar\alpha\theta\varphi_1)&=c,
%\\P^+(\bar\alpha\varphi_2- {\overline{\alpha(0)}}\varphi_4)&= c,
%\endaligned
%\right.
%\end{equation}
%which can also be written as
\begin{equation}\label{c18}
\left\{
\aligned
P^-(\varphi_4)&={\overline{\theta(0)}}P^-(\theta\varphi_3),
\\
P^-(\bar\alpha\theta\varphi_1)&={\overline{\alpha(0)}}P^-(\theta\varphi_3),\\
 P^+(\bar\alpha\varphi_2)&={\overline{\theta(0)}}P^+(\bar\alpha\theta\varphi_1)+c,
\\P^+(\bar\alpha\varphi_2)&= {\overline{\alpha(0)}}P^+(\varphi_4)+ c,
\endaligned
\right.\quad c\in\mathbb{C}.
\end{equation}
We now consider three possible cases.

{\bf  Case (i).} If ${{\theta(0)}}\ne 0$, then \eqref{c18} is equivalent to
\begin{equation}\label{c19}
\left\{
\aligned
P^+(\bar\alpha\theta\varphi_1)&=\tfrac{\overline{\alpha(0)}}{\overline{\theta(0)}}P^+(\varphi_4),\\
P^-(\bar\alpha\theta\varphi_1)&=\tfrac{\overline{\alpha(0)}}{\overline{\theta(0)}}P^-(\varphi_4),\\
 P^+(\bar\alpha\varphi_2)&={\overline{\alpha(0)}}P^+(\varphi_4)+c,
\\P^-(\theta\varphi_3)&=\tfrac{1}{\overline{\theta(0)}}P^-(\varphi_4).
\endaligned
\right.
\end{equation}
In other words,
\begin{equation}\label{c20}
\left\{
\aligned
\varphi_1&=\tfrac{\overline{\alpha(0)}}{\overline{\theta(0)}}\alpha\bar\theta\varphi_4,\\
 \varphi_2&={\overline{\alpha(0)}}\alpha P^+(\varphi_4)+c\alpha+\alpha P^-(\bar\alpha \varphi_2),
\\\varphi_3&=\bar \theta P^+(\theta\varphi_3)+\tfrac{1}{\overline{\theta(0)}}\bar \theta P^-(\varphi_4)
\endaligned
\right.
\end{equation}
and $\varphi_4\in L^{\infty}$ is arbitrary. Since
$$c\alpha+\alpha P^-(\bar\alpha \varphi_2)\perp \alpha z H^2\quad\text{and}\quad \bar \theta P^+(\theta\varphi_3)\perp \overline{\theta z H^2},$$
%the result follows.
we see that
$$D=\begin{bmatrix}\hat{T}_{\varphi_1}^{\theta,\alpha}&\check{\Gamma}_{\varphi_2}^{\alpha}
\\\hat{\Gamma}_{\varphi_3}^{\theta}&\check{T}_{\varphi_4}\end{bmatrix}=\begin{bmatrix}\hat{T}_{\psi_1}^{\theta,\alpha}&\check{\Gamma}_{\psi_2}^{\alpha}
\\\hat{\Gamma}_{\psi_3}^{\theta}&\check{T}_{\psi_4}\end{bmatrix}$$
with arbitrary $\psi_4=\varphi_4\in L^\infty$ and
\begin{equation*}
\left\{
\aligned
\psi_1&=\tfrac{\overline{\alpha(0)}}{\overline{\theta(0)}}\alpha\bar\theta\psi_4,\\
\psi_2&={\overline{\alpha(0)}}\alpha P^+(\psi_4),
\\\psi_3&=\tfrac{1}{\overline{\theta(0)}}\bar \theta P^-(\psi_4).
\endaligned
\right.
\end{equation*}
Note that although $\varphi_2$ and $\varphi_3$ are bounded the same is not necessarily true for $\psi_2$ and $\psi_3$.

{\bf  Case (ii).} If ${{\theta(0)}}=0$ and ${{\alpha(0)}}\ne 0$, then \eqref{c18} is equivalent to
\begin{equation}\label{c21}
\left\{
\aligned
P^-(\theta\varphi_3)&=\tfrac{1}{\overline{\alpha(0)}}P^-(\bar\alpha\theta\varphi_1),\\
P^-(\varphi_4)&=P^+(\varphi_4)=0,\\
 P^+(\bar\alpha\varphi_2)&=c,
\endaligned
\right.
\end{equation}
that is,
\begin{equation}\label{c22}
\left\{
\aligned
 \varphi_2&=c\alpha+\alpha P^-(\bar\alpha \varphi_2),
\\\varphi_3&=\bar \theta P^+(\theta\varphi_3)+\tfrac{1}{\overline{\alpha(0)}}\bar \theta P^-(\bar\alpha \theta \varphi_1),\\
\varphi_4&=0
\endaligned
\right.
\end{equation}
and $\varphi_1\in L^{\infty}$ is arbitrary. Here
$$\varphi_2\perp \alpha z H^2\quad\text{and}\quad \bar \theta P^+(\theta\varphi_3)\perp \overline{\theta z H^2},$$
so we can take arbitrary $\psi_1\in L^\infty$, $\psi_2=\psi_4=0$ and $\psi_3=\tfrac{1}{\overline{\alpha(0)}}\bar \theta P^-(\bar\alpha \theta \psi_1)$.

{\bf  Case (iii).} If ${{\theta(0)}}={{\alpha(0)}}= 0$, then \eqref{c18} is equivalent to
\begin{equation}\label{c23}
\left\{
\aligned
\bar\alpha\theta\varphi_1&\in H^{\infty},\\
\varphi_4&\in H^{\infty},\\
 \bar\alpha\varphi_2&\in\overline{zH^{\infty}}+\mathbb{C},
\endaligned
\right.
\end{equation}
with $\varphi_3\in L^{\infty}$ arbitrary. In this case
$$\varphi_2\perp \alpha z H^2$$
and we can take arbitrary $\psi_1\in \bar\theta \alpha H^\infty$, $\psi_4\in H^\infty$, $\psi_3\in \overline{\theta z H^2}$ and $\psi_2=0$.
\end{proof}

\begin{example}
	Let $\theta$ and $\alpha$ be two arbitrary nonconstant inner functions.
	\begin{enumerate}[(a)]
		\item If $\theta(0)\neq0$ and $\alpha(0)\ne 0$ (for simplicity of notations assume also that $\alpha(0)=\theta(0)$), then we can for example take $\psi_4=\theta$ or $\psi_4=\bar z$ or $\psi_4=\bar z \theta$ in \eqref{c1c2n}. Then we obtain the following operators from $\mathcal{I}(K_\theta^\perp,K_\alpha^\perp)$:
		\begin{equation}\label{d1a} \begin{bmatrix}\hat{T}_{\alpha}^{\theta,\alpha}&\check{\Gamma}_{\overline{\theta(0)}\alpha\theta}^{\alpha}
		\\0&\check{T}_{\theta}\end{bmatrix},\
 \begin{bmatrix}\hat{T}_{\bar z\alpha\bar \theta}^{\theta,\alpha}&0\\
\hat{\Gamma}_{\frac{1}{ \overline{\theta(0)}}\bar z\bar\theta}^{\theta}
		&\check{T}_{\bar z}\end{bmatrix},\
\begin{bmatrix}\hat{T}_{\bar z\alpha}^{\theta,\alpha}&\check{\Gamma}_{\overline{\theta(0)}\bar z\bar\alpha(\theta-\theta(0))}^{\alpha}
	\\
\hat{\Gamma}_{\frac{\theta(0)}{ \overline{\theta(0)}}\bar z\bar\theta}^{\theta}
		&\check{T}_{\bar z\theta}\end{bmatrix}.
		\end{equation}%
%		\begin{equation}\label{d1a} \begin{bmatrix}\hat{T}_{\frac{\overline{\alpha(0)}}{\overline{\theta(0)}}\alpha}^{\theta,\alpha}&\check{\Gamma}_{\overline{\alpha(0)}\alpha\theta}^{\alpha}
%		\\0&\check{T}_{\theta}\end{bmatrix},
% \begin{bmatrix}\hat{T}_{\frac{\overline{\alpha(0)}}{\overline{\theta(0)}}\bar z\alpha\bar \theta}^{\theta,\alpha}&0\\
%\hat{\Gamma}_{\frac{1}{ \overline{\theta(0)}}\bar z\alpha\bar\theta}^{\theta}
%		&\check{T}_{\bar z}\end{bmatrix},
%\begin{bmatrix}\hat{T}_{\frac{\overline{\alpha(0)}}{\overline{\theta(0)}}\bar z\alpha}^{\theta,\alpha}&\check{\Gamma}_{\overline{\alpha(0)}\bar z\bar\alpha(\theta-\theta(0))}^{\alpha}
%	\\
%\hat{\Gamma}_{\frac{\theta(0)}{ \overline{\theta(0)}}\bar z\bar\theta}^{\theta}
%		&\check{T}_{\bar z\theta}\end{bmatrix}
%		\end{equation}

\item If  $\theta(0)\neq0$ and $\alpha(0)=0$, then choosing $\psi_4$ as in (a) we obtain operators
		\begin{equation}\label{d1b} \begin{bmatrix}0&0
		\\0&\check{T}_{\theta}\end{bmatrix},\quad
 \begin{bmatrix}0&0\\
\hat{\Gamma}_{\frac{1}{ \overline{\theta(0)}}\bar z\bar\theta}^{\theta}
		&\check{T}_{\bar z}\end{bmatrix},\quad
\begin{bmatrix}0&0
	\\
\hat{\Gamma}_{\frac{\theta(0)}{ \overline{\theta(0)}}\bar z\bar\theta}^{\theta}
		&\check{T}_{\bar z\theta}\end{bmatrix}
		\end{equation}
\item If $\theta(0)=0$ and $\alpha(0)\ne 0$, then we can for example take $\psi_1=\alpha$ or $\psi_1=\bar \theta$ or $\psi_1=\bar z \bar \theta$ in \eqref{fi32}. Then we obtain the following operators
		\begin{equation}\label{d2}
		\begin{bmatrix}\hat{T}_{\alpha}^{\theta,\alpha}&0
		\\0&0\end{bmatrix},\quad \begin{bmatrix}\hat{T}_{\bar\theta}^{\theta,\alpha}&0
		\\\hat\Gamma^{\theta}_{\bar\theta(\overline{\frac{\alpha}{\alpha(0)}}-1)}&0\end{bmatrix},
\quad \begin{bmatrix}\hat{T}_{\bar z\bar\theta}^{\theta,\alpha}&0
		\\\hat\Gamma^{\theta}_{\frac{1}{\overline{\alpha(0)}}\bar z \bar \alpha\bar\theta}&0\end{bmatrix}.
		\end{equation}
		
			\item If $\theta(0)=0=\alpha(0)$
 and $\chi_1, \chi_3,\chi_4\in H^\infty$ are arbitrary, then any linear combination of
	\begin{equation}\label{de4}\begin{bmatrix}\hat{T}_{\alpha\bar\theta\chi_1}^{\theta,\alpha}&0
\\0&0\end{bmatrix}, \quad\begin{bmatrix}0&0
\\\hat{\Gamma}_{\bar \theta \bar\chi_3}^{\theta}&0\end{bmatrix},\quad\begin{bmatrix}0&0
\\0&\check{T}_{\chi_4}\end{bmatrix}
		\end{equation}
satisfies \eqref{inttwin} and belongs to $\mathcal{I}(K_\theta^\perp,K_\alpha^\perp)$.

	\end{enumerate}
	\end{example}

We can now use Theorem \ref{kmutant} to describe all those operators $D\in \mathcal{B}(K_\theta^\perp,K_\alpha^\perp)$ that satisfy \eqref{inttwin} and belong to $\mathcal{T}(K_\theta^\perp,K_\alpha^\perp)$ (give another proof of Theorem \ref{idatto}). Clearly, $D=0$ is one such operator. We can assume that $D\neq 0$.

By Theorem \ref{kmutant}, $D$ satisfies \eqref{inttwin} if and only if
\begin{equation}\label{raz}
D=\begin{bmatrix}\hat{T}_{\psi_1}^{\theta,\alpha}&\check{\Gamma}_{\psi_2}^{\alpha}
\\\hat{\Gamma}_{\psi_3}^{\theta}&\check{T}_{\psi_4}\end{bmatrix}
\end{equation}
with $\psi_i$, $i=1,2,3,4$, satisfying conditions \eqref{c1c2n}, \eqref{fi32} or \eqref{fi33} (depending on the case). On the other hand, $D\in \mathcal{T}(K_\theta^\perp,K_\alpha^\perp)$, $D=D_{\varphi}^{\theta,\alpha}$ for some $\varphi\in L^{\infty}$, if and only if
\begin{equation}\label{dwa}
D=\begin{bmatrix}\hat{T}_{\varphi}^{\theta,\alpha}&\check{\Gamma}_{\varphi}^{\alpha}
\\\hat{\Gamma}_{\varphi}^{\theta}&\check{T}_{\varphi}\end{bmatrix}.
\end{equation}
Comparing \eqref{raz} and \eqref{dwa} we get
\begin{equation}\label{trzy}
\hat{T}_{\varphi}^{\theta,\alpha}=\hat{T}_{\psi_1}^{\theta,\alpha},\quad \check{T}_{\varphi}=\check{T}_{\psi_4},\quad \check{\Gamma}_{\varphi}^{\alpha}=\check{\Gamma}_{\psi_2}^{\alpha}\quad\text{and}\quad \hat{\Gamma}_{\varphi}^{\theta}=\hat{\Gamma}_{\psi_3}^{\theta}.
\end{equation}
First two equalities in \eqref{trzy} amount to
\begin{equation}\label{cztery}
\varphi=\psi_1=\psi_4,
\end{equation}
while the last two give
\begin{equation}\label{piec}
P^-(\alpha\bar \varphi)=P^-(\alpha\bar \psi_2)
%\end{equation}
\qquad\text{and}\qquad
%\begin{equation}\label{szesc}
P^-(\theta \varphi)=P^-(\theta \psi_3).
\end{equation}
By \eqref{cztery}--\eqref{piec}, from the cases considered in Theorem \ref{kmutant} only cases (i) and (iii) are possible ((ii) is not possible since $D\neq 0$).

{\bf  Case (i).} Here ${{\theta(0)}}\ne 0$, $\psi_4\in L^{\infty}$ and
\begin{equation}\label{c29}
\left\{
\aligned
\psi_1&=\tfrac{\overline{\alpha(0)}}{\overline{\theta(0)}}\alpha\bar\theta\psi_4,\\
 \psi_2&={\overline{\alpha(0)}}\alpha P^+(\psi_4),
\\\psi_3&=\tfrac{1}{\overline{\theta(0)}}\bar \theta P^-(\psi_4).
\endaligned
\right.
\end{equation}
The first of the three equalities above together with \eqref{cztery} gives
$$1=\tfrac{\overline{\alpha(0)}}{\overline{\theta(0)}}\alpha\bar\theta.$$
Clearly, then $\alpha(0)\neq 0$. It follows that $\alpha$ is a constant multiple of $\theta$, $\alpha=\lambda\cdot\theta$, $\lambda\in\mathbb{T}$.

By \eqref{piec} and \eqref{c29},
$$P^-(\theta \varphi)=P^-(\theta \psi_3)=\tfrac{1}{\overline{\theta(0)}}  P^-(\psi_4)=\tfrac{1}{\overline{\theta(0)}}  P^-(\varphi),$$
which means that
$$\varphi-\overline{\theta(0)}\theta\varphi=k_0^{\theta}\varphi \in H^2.$$
Thus $\varphi\in L^{\infty}\cap H^2=H^{\infty}$ and there exists $h\in H^{\infty}$ such that $\theta\varphi=h$, $\varphi=h/\theta=\bar\theta h$. Since $\alpha=\lambda\cdot\theta$ and $\varphi=\psi_4\in H^\infty$, by \eqref{piec} and \eqref{c29} we get
$$P^-(\theta^2 \bar h)=P^-(\theta(0)\theta \bar h).$$
Hence, for each $f\in H^2$,
\begin{align*}
\langle k_0^{\theta}\varphi,\theta z f\rangle &=\langle \varphi-\overline{\theta(0)}\theta\varphi,\theta z f\rangle =\langle \bar\theta h,\theta z f\rangle -\langle \overline{\theta(0)}\bar\theta h, z f\rangle\\
&=\langle \overline{z f},\theta^2 \bar h\rangle -\langle \overline{z f} , {\theta(0)}\theta\bar h\rangle=0
\end{align*}
and $k_0^{\theta}\varphi\in K_{z\theta}$.

{\bf  Case (iii).} Here ${{\theta(0)}}={{\alpha(0)}}= 0$ and
\begin{equation}%\label{c23}
\psi_1\in \alpha\bar\theta H^{\infty},\quad \psi_2=0,\quad \psi_3\in L^{\infty}\quad\text{and}\quad \psi_4\in H^{\infty}.
\end{equation}
Denote $\theta=\theta_1\cdot \textrm{gcd}(\alpha,\theta)$. Since there exists $h\in H^{\infty}$ such that
$$\varphi=\alpha \bar\theta h=\psi_4\in H^{\infty},$$
we must have
$$\varphi=\tfrac{\alpha}{\textrm{gcd}(\alpha,\theta)} \cdot\tfrac{h}{\theta_1}\in \tfrac{\alpha}{\textrm{gcd}(\alpha,\theta)} H^{\infty}.$$
Moreover, by \eqref{piec},
$$P^-(\theta\bar h)=0$$
and so, for each $f\in H^2$,
$$\langle h/\theta_1,\textrm{gcd}(\alpha,\theta)\cdot zf\rangle =\langle h,\theta zf\rangle =\langle \overline{zf}, \theta\bar h\rangle =0.$$
That is, $h/\theta_1\in K_{z\cdot \textrm{gcd}(\alpha,\theta)}$ and
$$\varphi\in \tfrac{\alpha}{\textrm{gcd}(\alpha,\theta)} K_{z\cdot \textrm{gcd}(\alpha,\theta)}.$$

Note that by Theorem \ref{idatto} the set $\mathcal{I}(K_\theta^\perp,K_\alpha^\perp)$ contains operators from $\mathcal{T}(K_\theta^\perp,K_\alpha^\perp)$ if and only if $\alpha=\lambda\cdot\theta$ or $\alpha(0)=\theta(0)=0$. Otherwise no operators from $\mathcal{T}(K_\theta^\perp,K_\alpha^\perp)$ satisfy \eqref{inttwin}.

\begin{example}
Take any nonconstant inner function $\alpha$ and let $\theta =z\alpha$. Then $\textrm{gcd}(\alpha,\theta)=\alpha$ (here $\alpha\leqslant \theta$) and
$$\tfrac{\alpha}{\textrm{gcd}(\alpha,\theta)} K_{z\cdot \textrm{gcd}(\alpha,\theta)}=K_{\theta}=K_{z\alpha}=K_{z}\oplus zK_{\alpha}.$$  If $\alpha(0)=0$, then by Theorem \ref{idatto} every $D_{\varphi}^{\theta,\alpha}$ with $\varphi\in K_{\theta}\cap L^{\infty}$ satisfies \eqref{inttwin} (note that $K_{\theta}\cap L^{\infty}$ is always nonempty; here it contains for example constant functions). On the other hand, if $\alpha(0)\neq 0$, then no operator from $\mathcal{T}(K_\theta^\perp,K_\alpha^\perp)$ satisfies \eqref{inttwin}.
\end{example}

\begin{remark}Can we find nonconstant $\alpha$ and $\theta$ such that every operator satisfying \eqref{inttwin} belongs to $\mathcal{T}(K_\theta^\perp,K_\alpha^\perp)$? The answer is no. If $\alpha(0)=\theta(0)=0$, then for example the first operator in \eqref{d1b} satisfies \eqref{inttwin} but does not belong to $\mathcal{T}(K_\theta^\perp,K_\alpha^\perp)$ (as $\theta\neq 0$). If $\alpha=\lambda \cdot\theta$ and $\theta(0)\neq 0$ (otherwise $\alpha(0)=\theta(0)=0$) we can consider the first operator in \eqref{d1a}, which in this case takes the form
	\begin{equation*}%\label{d1}
D=\begin{bmatrix}\hat{T}_{\theta}^{\theta,\alpha}&\check{\Gamma}_{\overline{\theta(0)}\theta^2}^{\alpha}
\\0&\check{T}_{\theta}\end{bmatrix}.
\end{equation*}
Again, $D$ satisfies \eqref{inttwin} but does not belong to $\mathcal{T}(K_\theta^\perp,K_\alpha^\perp)$. If it did belong to $\mathcal{T}(K_\theta^\perp,K_\alpha^\perp)$, then we would have $D=D^{\theta,\alpha}_{\theta}$ and by Theorem \ref{idatto}, $\theta k_0^{\theta}\in K_{z\theta}$. But since for each $h\in H^2$,
$$\langle \theta k_0^{\theta} ,\theta zh\rangle=\langle 1 ,zh\rangle-\overline{\theta(0)}\langle \theta ,zh\rangle=-\overline{\theta(0)}\langle \theta ,zh\rangle$$
and $\theta(0)\neq 0$, we see that $\theta k_0^{\theta}\in K_{z\theta}$ if and only if $\theta$ is a constant function which is a contradiction.
\end{remark}

\end{document}